\documentclass[12pt]{article}           
\long\def\new#1\endnew{{\bf #1}}    
\long\def\del#1\enddel{} 
\def\preprint{HUB-EP-00/03\\TUW--00/01}
\def\finished{January 2000}
\def\archive {math.AG/0001106}           

\def\Title{Reflexive polyhedra, weights and\\[4mm]toric Calabi--Yau fibrations}

\long\def\Abstract{
During the last years we have generated a large number of data related to 
Calabi--Yau hypersurfaces in toric varieties which can be described by 
reflexive polyhedra. We classified all reflexive polyhedra in three dimensions 
leading to K3 hypersurfaces and have nearly completed the four dimensional case
relevant to Calabi--Yau threefolds.
In addition, we have analysed for many of the resulting spaces whether they
allow fibration structures of the types that are relevant in the context of
superstring dualities.
In this survey we want to give background information both on how we obtained
these data, which can be found at our web site, and on how they may be used.
We give a complete exposition of our classification algorithm at a
mathematical (rather than algorithmic) level. We also describe how fibration
structures manifest themselves in terms of toric diagrams and how we 
managed to find the respective data.
Both for our classification scheme and for simple descriptions of fibration
structures the concept of weight systems plays an important role.
}
\def\cont#1{\mathop{\vtop{\ialign{##\crcr $\hfil\displaystyle
                     {#1}\hfil$\crcr\noalign{\kern3 pt \nointerlineskip}
                     $\bracelu\leaders\vrule\hfill\leaders\vrule\hfill
                     \braceru$\crcr\noalign{\kern3 pt }}}}\limits}

\def\CY{Calabi--Yau}

\def\ipo{\hbox{\bf 0}} 
\def\fib{_{\rm f}} 
\def\bas{_{\rm b}} 
\def\Kt{{\tilde K}}

\def\Def{\ni{\bf Definition: }}
\def\Pro{{\it Proof: }} 
\def\hBo{\hfill$\Box$}
\def\q{\hbox{\bf q}}
\def\ipo{\hbox{\bf 0}}

\usepackage{cite,epic,epsf,amssymb,latexsym}

\def\ifundefined#1{\expandafter\ifx\csname#1\endcsname\relax}
\def\bye{\end{document}}        

\def\HS#1 {\hspace*{#1pt}} \def\VS#1 {\vspace*{#1pt}} 
\def\VR#1#2{\vrule height #1mm depth #2mm width 0pt}     
\def\TVR#1#2{@{~~\VR{#1}{#2}}}                          
\def\BC{\begin{center}} \def\EC{\end{center}}

\let\0=\over      \let\1=\vec      \def\2{{1\over2}}
\let\4=\underline \let\5=\overline \let\6=\partial      \def\7#1{{#1}\llap{/}}
\def\8#1{{\textstyle{#1}}}         \def\9#1{{\ifmmode{\pmb{#1}}\else\bf#1\fi}}

      \def\({\left(}       \def\){\right)}

\def\eeql#1 {\label{#1}\eeq}      \let\nn=\nonumber  
\def\beq{\begin{equation}}      \def\eeq{\end{equation}}        
\def\bea{\begin{eqnarray}}      \def\eea{\end{eqnarray}} 

\def\mao#1{\mathop{\rm #1}\nolimits}  
       \def\Tr{\mao{Tr}}

\let\and=\wedge

\let\bra=\langle        \let\ket=\rangle        \def\<#1\>{\bra #1 \ket}
\let\ni=\noindent

\def\rel#1 #2{\buildrel #1 \over {#2}}  
\def\fnote#1#2{\begingroup\def\thefootnote{#1}\footnote{#2}
                \addtocounter{footnote}{-1}\endgroup}   

\def\subdef#1{\gdef\globalColor##1{##1}}      
      \subdef{Black}

          \let\d=\delta   
         \let\th=\theta  \let\e=\varepsilon
   \let\l=\lambda        
            \let\p=\pi            
                 
               \let\S=\Sigma 
\let\P=\Pi             \let\D=\Delta

\def\ce{{\cal E}}

\def\IR{{\mathbb R}} \def\IC{{\mathbb C}} \def\IP{{\mathbb P}} 
  \def\IN{{\mathbb N}} 
\def\IZ{{\mathbb Z}}

\def\plb#1 #2 {Phys. Lett. {\bf B#1} #2 }
\def\phr#1 #2 {Phys. Rep. {\bf  #1} #2 }        
\def\npb#1 #2 {Nucl. Phys. {\bf B#1} #2 }
\def\aph#1 #2 {Ann. Phys. {\bf #1} #2 }         
\def\jmp#1 #2 {J. Math. Phys. {\bf #1} #2 }
\def\jgp#1 #2 {J. Geom. Phys. {\bf #1} #2 }
\def\prd#1 #2 {Phys. Rev. {\bf D#1} #2 }
\def\prl#1 #2 {Phys. Rev. Lett. {\bf #1} #2 }
\def\rmp#1 #2 {Rev. Mod. Phys.  {\bf #1} #2 }
\def\zpc#1 {Z. Phys. {\bf #1C} }
\def\cmp#1 #2 {Commun. Math. Phys. {\bf #1} #2 }
\def\cqg#1 #2 {Class.Quant.Grav. {\bf #1} #2 }
\def\mpl#1 {Mod. Phys. Lett. {\bf A#1} }
\def\cpc#1 {Computer Phys. Commun. {\bf #1} }   
\def\ijmp#1 {Int. J. Mod. Phys. {\bf A#1} }
\def\ijmpC#1 {Int. J. Mod. Phys. {\bf C#1} }

\def\BP{\begin{picture}} \def\EP{\end{picture}}         
\newcounter{TRefNX} \let\OLDcite=\cite  \makeatletter
\def\makeTRefs#1{\@for  \NewTRef:=#1\do{\global\makeTRef{\NewTRef}}}
\def\makeTRef#1{\ifundefined{TRef#1}\stepcounter{TRefNX}%
\expandafter\xdef\csname TRef#1\endcsname{\theTRefNX}\fi}\makeatother
\def\NEWcite#1{\makeTRefs{#1}\OLDcite{#1}}  

\ifundefined{draftmode} {} \else        \let\cite=\NEWcite
   \newcount\HOUR  \HOUR=\the\time \divide\HOUR by 60    \multiply \HOUR by 60 
   \newcount\MIN   \MIN=\the\time  \advance\MIN by -\HOUR   \divide\HOUR by 60
   \ifnum   \MIN>9  \def\printTIME{{\it\the\HOUR\,:\,\the\MIN}}
         \else   \def\printTIME{{\it\the\HOUR\,:\,0\the\MIN}} \fi 
   
   \def\LLab#1{\BP(0,0)\unitlength=1mm\put(-12,.5){\makebox(0,0)[cr]{\small #1
        \rlap{$_{_{\makeatletter\csname TRef#1\endcsname\makeatother}}$}}}\EP}
\fi

\begin{document}
{
\vspace*{-78pt}\begin{flushright}  \archive\\ \preprint  \end{flushright}
\vspace{-8mm}
\BC{\huge\bf \Title}
\\[5mm]
        Maximilian KREUZER\fnote{*}{e-mail: kreuzer@hep.itp.tuwien.ac.at}
\\[2mm]
        Institut f\"ur Theoretische Physik, Technische Universit\"at Wien\\
        Wiedner Hauptstra\ss e 8--10, A-1040 Wien, AUSTRIA
\\[2mm]                       and
\\[2mm] Harald SKARKE\fnote{\#}{e-mail: skarke@physik.hu-berlin.de}
\\[1mm] Humboldt Universit\"at zu Berlin, Institut f\"ur Physik, QFT,
\\      Invalidenstra\ss e 110, D-10115 Berlin, GERMANY
\\[5mm]       
{\bf ABSTRACT } \\[4mm]  \parbox{15cm}
{\baselineskip=14.5pt ~~~\Abstract} \EC
\vfil \noindent \preprint \\ \finished 
\setcounter{page}{0} \thispagestyle{empty} \newpage \pagestyle{plain}  
}

\section{Introduction}          

A few years after Calabi--Yau manifolds had found their way into physics
it was conjectured 
that they should actually come in pairs with opposite
Euler number, since an exchange of complex structure and K\"ahler moduli 
in physics corresponds to a change of sign in the definition 
of the charge, or, equivalently, an exchange of particles and
anti-particles \cite{dix,lvw}. This phenomenon is called mirror symmetry.
Although the situation is 
complicated by the fact that 
there are
rigid Calabi--Yau manifolds whose ``mirror string compactifications'' do 
not have a straightforward geometrical interpretation \cite{rigid,AG}, 
the search for the 
mirror manifolds proved to be an extremely fruitful enterprise from both the 
physicists' and the  mathematicians' perspective \cite{exact1,exact2,lagr}.

The first systematic constructions of large classes of Calabi--Yau 
threefolds as complete intersections in products of projective spaces 
\cite{CICY} did not seem to support the mirror hypothesis 
because the resulting manifolds all had negative Euler numbers. 
But when the attention was extended to weighted projective spaces, it turned
out that the blow up parameters
of the quotient singularities can provide large positive contributions. 
The first substantial list of pairs of Hodge numbers resulting from
constructions of this type \cite{CLS} was almost mirror symmetric in the 
sense that only for a few percent of the Hodge data
the respective mirror pair was not in the list.
A complete classification \cite{nms,KlS}, however,
made the picture worse, and  abelian quotients \cite{aas}, which make a 
subclass of these spaces perfectly symmetric \cite{BH,odt}, 
did not help with this problem either. 

Batyrev's construction of toric Calabi--Yau hypersurfaces
\cite{bat}, which is manifestly mirror symmetric while generalizing the 
above results, provided a solution to this puzzle. 
In this framework the geometrical
data is encoded by a reflexive polyhedron, i.e. a lattice polyhedron whose
facets are all at distance 1 from the unique interior point (see below).
Toric geometry turned out to provide a very 
efficient tool for the analysis of many physical aspects of Calabi--Yau 
compactifications, including the physics of perturbative \cite{agm}
and non-perturbative \cite{cgh1,cgh2,Str,gms} topology changing transitions, 
as well
as fibration structures that are important in string dualities \cite{dual,Fv}.

This made a constructive classification of reflexive polyhedra
a useful and interesting enterprise. Our 
approach to this problem \cite{crp} was partly inspired by our experience 
with the classification of weighted projective spaces 
that admit transversal quasi-homogeneous polynomials \cite{cqf}. 
Indeed, as it turned out, the Newton polyhedra that correspond
to polynomials defining CY hypersurfaces in weighted $\IP^4$ are all reflexive
\cite{cok,wtc} and provide a canonical resolution of the ambient space 
singularities
(this is no longer true in higher dimensions). Actually, regardless of
the transversality condition, a diophantine equation of the form 
$\sum n_ia_i=d$ with positive coefficients $n_i$, $d=\sum n_i$, and with the 
set of solutions restricted to $a_i\ge0$ gives a simple way to produce lattice
polyhedra with at most one interior point (this is a necessary condition 
for reflexivity): We may regard this as an embedding of the lattice into
a higher dimensional space with the polyhedron being contained in the finite
intersection of an affine subspace with the non-negative half-spaces. 
All lattice points, except for the candidate interior point, whose 
coordinates are all equal to 1, 
are located on some coordinate hyperplane $a_i=0$.

We may then ask ourselves if all reflexive polyhedra are contained in 
polyhedra that can be embedded in this way. 
In the next section we will show that the answer is assertive provided 
that we allow for an embedding with higher codimension $k-n$, i.e. we also 
consider solutions to more than one equation of the above form, 
\beq
        \sum_{i=1}^k a_i n_i^{(j)} = d^{(j)},           \label{AffineEq}
                ~~~~~ d^{(j)}=\sum_{i=1}^kn_i^{(j)}, ~~~~~  j=1,\ldots,k-n,
\eeq
but with some of the coefficients $n_i^{(j)}$ equal to
zero according to a certain pattern. What then makes our
construction work is the fact that there is only a finite set of coefficients
that lead to lattice polytopes with an interior point. The collections 
$n_i^{(j)}$ of non-negative numbers are called weight systems in the case of 
a single equation and combined weight systems if $k-n>1$.
If we shift our coordinates to $x_i=a_i-1$ the resulting polyhedron lies
in a linear subspace of the embedding space determined by 
$\sum_i n_i^{(j)}x_i=0$, is bounded by $x_i\ge-1$ and has the origin of the 
embedding space as its interior point.
These linear coordinates are more useful for many
general considerations wheras the affine coordinates $a_i$ are better suited 
for quickly finding the lattice points in a given example.
\del
We will mostly work with shifted coordinates $x_i=a_i-1$ that move the 
interior point into the origin of the embedding space and, accordingly, with
linear equations $\sum_i n_i^{(j)}x_i=0$ and $x_i\ge-1$. 
\enddel

\begin{figure}
\BC     \newcount\XYfac \newcount\Yfac  \newcount\YA    \linethickness{3pt}
        \newcount\XZfac \newcount\Zfac  \newcount\ZA    \thicklines
        \newcount\Xaux  \newcount\Yaux  \newcount\Zaux
\def\makeYZ(#1,#2,#3){
        \Xaux=#1 \multiply \Xaux by \XYfac \Yaux=#2 \multiply \Yaux by \Yfac 
        \advance\Yaux by -\Xaux
        \Xaux=#1 \multiply \Xaux by \XZfac \Zaux=#3 \multiply \Zaux by \Zfac 
        \advance\Zaux by -\Xaux
}
\def\xyzline(#1,#2,#3)(#4,#5,#6){\makeYZ(#1,#2,#3) \YA=\Yaux \ZA=\Zaux 
        \makeYZ(#4,#5,#6)       \drawline(\YA,\ZA)(\Yaux,\Zaux)
}       
\def\xyzdash(#1,#2,#3)(#4,#5,#6){\makeYZ(#1,#2,#3) \YA=\Yaux \ZA=\Zaux 
        \makeYZ(#4,#5,#6)       \dashline2(\YA,\ZA)(\Yaux,\Zaux)
}
\def\xyzput(#1,#2,#3)#4{\makeYZ(#1,#2,#3)\put(\Yaux,\Zaux){#4}}
\def\xyzdot(#1,#2,#3){\xyzput(#1,#2,#3)\pt}
\unitlength=1mm         \unitlength=2pt \XYfac=7  \Yfac=20  \XZfac=5 \Zfac=30
                                        \def\pt{{\circle*3}}
\BP(120,80)(-60,-15) 
\put(-60,0){    \xyzdash(-1,0,0)(0,1,0)\xyzline(1,0,0)(0,-1,0)
                \xyzline(1,0,0)(0,1,0)\xyzdash(-1,0,0)(0,-1,0)
        \xyzline(1,0,0)(0,0,2)\xyzline(0,-1,0)(0,0,2)
        \xyzline(0,1,0)(0,0,2)\xyzdash(-1,0,0)(0,0,2)

        \makeYZ(-1,1,0) \YA=\Yaux \ZA=\Zaux \makeYZ(-1,0,0)
        \multiput(0,0)(0,\Zfac)3\pt
        \xyzdot(1,0,0)\xyzdot(-1,0,0)\xyzdot(0,1,0)\xyzdot(0,-1,0)
   \xyzput(0,-1,0){{\put(-6,-2)1}}\xyzput(0,1,0){{\put(3,-2)2}}
   \xyzput(1,0,0){{\put(-5,-5)3}}\xyzput(-1,0,0){{\put(2,2)4}}
   \xyzput(0,0,2){{\put(2.5,0)5}}\xyzput(0,-2,2){{\put(5,-15){$\nabla=\D^*$}}}
}
\put(50,0){     \xyzline(-2,2,0)(2,2,0)\xyzline(2,-2,0)(-2,-2,0)
        \xyzline(2,-2,0)(2,2,0)\xyzline(-2,2,0)(-2,-2,0)
        \xyzline(2,-2,0)(0,0,2)\xyzline(-2,2,0)(0,0,2)
        \xyzline(2,2,0)(0,0,2)\xyzdash(-2,-2,0)(0,0,2)
                                \xyzput(0,-2,2){{\put(5,-15){$\D$}}}    

        \multiput(0,0)(0,\Zfac)3\pt
        \makeYZ(-1,1,0) \YA=\Yaux \ZA=\Zaux \makeYZ(-1,0,0)
        \put(0,\Zfac){\matrixput(-\YA,-\ZA)(\Yfac,0)3(\Yaux,\Zaux)3\pt}
        \makeYZ(-2,2,0) \YA=\Yaux \ZA=\Zaux \makeYZ(-1,0,0)
        \matrixput(-\YA,-\ZA)(\Yfac,0)5(\Yaux,\Zaux)5\pt
}
\EP
\\[3pt]{\bf Fig. 1:} A minimal polyhedron $\nabla$ that corresponds to a 
        combined weight system.
\EC
\end{figure}

Let us illustrate with the example in {\it fig.} 1 how we can obtain a weight 
system for a given reflexive polyhedron $\D$ with vertices in some 
$n$-dimensional lattice $M$.
Reflexivity implies that the dual (or polar) polytope $\D^*$ defined in 
eq. (\ref{dual}) below has its vertices on the dual lattice 
$N=\mao{Hom}(M,\IZ)$. 
In our case $\D^*$ is already minimal in the sense
that we lose the interior point (IP) if we drop any of its vertices and take
the convex hull of the remaining vertices. 
The set of vertices of $\nabla=\D^*$ can be decomposed
into the two triangles $(V_1,V_2,V_5)$ and $(V_3,V_4,V_5)$ that both contain 
the IP in their lower dimensional interior. 
As we will show later, similar decompositions are always possible for 
minimal polyhedra.
For both triangles the 
barycentric coordinates of the IP are given by $\q=(1/4,1/4,1/2)$, i.e.
$\sum q_i V_i=0$ and $\sum q_i=1$, where the sum is over the indices of the
vertices for any of the two triangles. Rescaling the coefficients 
to integers $n_i^{(j)}=d^{(j)}q_i^{(j)}$ we arrive at the weight system 
$n_i^{(1)}=(1,1,0,0,2)$, $n_i^{(2)}=(0,0,1,1,2)$.
We will demonstrate in the next section that weights obtained in this way
can always be used to describe the dual polytope $\D$ as in eq. \ref{AffineEq}.
In the present case, this construction leads to
\bea    x_1+x_2~~~~~~~~~~~~+2x_5&=&0,\\
        ~~~~~~~~~~~x_3+x_4+2x_5&=&0.
\eea
Eliminating, for example, $x_2$ and $x_4$ it is easily checked that we indeed
reconstructed $\D$ (note that $a_i=x_i+1$ is the lattice distance of a point
from the facet dual to $V_i$).
If we keep all points with $x_i\ge-1$ then, in our example, $\D^*$ 
is equal to $\nabla$. In general $\D^*$ will not be minimal and we first have
to drop some vertices of $\D^*$ to arrive at a minimal polytope $\nabla$
whose simplex decomposition leads to a weight system.
If we drop points from $\D$ in such a way that $\D'\subset\D$ is reflexive, 
then $\D'^*$ becomes larger. The vertices of $\nabla$ remain vertices of 
$\D'^*\supset\nabla$ as long as the bounding hyperplanes $x_i=-1$, 
which in our case support all facets of $\D$, 
are affinely spanned by facets of $\D'$. 

A different way to generate a `smaller' $\D$ is to keep the vertices but to go 
to a coarser $M$ lattice:
We may, for example, demand $x_1-x_3\in 2\IZ$ or $x_1+x_5\in 2\IZ$.
Correspondingly, the $N$ lattice becomes finer and is no longer generated by
the vertices of $\nabla$. In general there will occur additional lattice 
points in $\nabla$.
The coarsest lattice that keeps all vertices of $\D$ and the IP is obtained by
imposing $x_1-x_3\in 4\IZ$ and $x_1+x_5\in 2\IZ$.
Actually, in our example, this exchanges $\nabla$ and $\D$.

In practice, because of the huge number of solutions, an enumeration of
all reflexive polyhedra seems to be possible only in up to 4 dimensions.
This leads to a further simplification of the procedure because in up to 
4 dimensions all polytopes $\D$ that correspond to a minimal $\nabla$ are 
reflexive \cite{wtc}. 
Moreover, $\D$ is contained in a larger polytope $\hat\D$ if
and only if $\hat\D^*$ is contained in $\D^*$. 
Therefore only 
minimal polytopes for which $\D^*$ does not contain any reflexive subpolytope
are necessary ingredients for our classification scheme.
We will show that in 4 dimension 
there are 308 
reflexive polytopes that contain all others as subpolytopes, provided that 
we also consider sublattices.
Finding all relevant lattices is a subtle point 
and our strategy to solve this problem will be described below. 
There are at least 25 additional
maximal reflexive polytopes that can be obtained from these 308 objects 
on sublattices.

While one of the main insights of the `first superstring revolution' was the
fact that Calabi--Yau spaces are crucial for string
compactifications, it was found during the `second string revolution' that 
fibration structures of Calabi--Yau manifolds are essential for understanding
various non-perturbative string dualities.
In particular, K3 fibrations are required for the duality between heterotic
and IIA theories \cite{KLM,AL} and elliptic fibrations are needed for F-theory
compactifications \cite{Fv,MVI,MVII}.
Again toric geometry provides beautiful tools for studying the 
respective structures.
As we will see, the polytope $\D^*\fib$ corresponding to the fiber manifests 
itself as a subpolytope of $\D^*$ with the same interior point, whereas
the base space is a toric variety whose fan can be determined by projecting
the original fan along the linear subspace spanned by $\D^*\fib$.
While we never attempted to give a complete classification of structures of
this type, we did create large lists of fibration structures \cite{k3,fft}.

Our data are accessible at our web site \cite{KScy}, and  
we plan to make the source code of our programs available in the near future.
Since one of the motivations for writing this contribution was to give 
useful background material for anyone interested in applying our data,
we would like to briefly mention some older results on our web page that 
will not be discussed in the remainder of this paper.
These are mostly 
related to weighted projective spaces and, in the physical context, to 
Landau-Ginzburg models \cite{lvw,va89}. We classified all 10839 weight systems 
allowing transversal quasihomogeneous polynomials \cite{cqf,nms} with 
singularity index 3, leading to Landau-Ginzburg models with a central charge
of c=9 and computed the corresponding numbers of (anti) chiral states in
the superconformal field theories (this includes the 7555 transversal weights
for weighted $\IP^4$).  Vafa's formulas for these numbers 
\cite{va89} inspired the definition of what Batyrev et al. call string
theoretic Hodge numbers \cite{BD}. We also extended these results to 
arbitrary abelian quotients that leave a transversal polynomial invariant 
\cite{aas} (and included the 
modifications by discrete torsions \cite{lgt}, 
which correspond to topologically non-trivial background 2-form fields in 
the physical context \cite{dt}). Since the Newton polyhedra are reflexive also 
for abelian quotients, the resulting Hodge numbers (without discrete torsion)
are all recovered in the toric context. Nevertheless our results might be 
useful when working in weighted projective spaces, since transversal 
polynomials in general have larger symmetries than the complete Newton 
polyhedra. 

We will not discuss Calabi--Yau data obtained by other groups here.
An important class of spaces that we did not consider consists of complete 
intersection Calabi--Yau varieties. 
The classification of these objects in products of projective spaces was given 
in \cite{CICY}, and Klemm has produced a sizeable list of codimension two 
complete intersections in weighted projective spaces which is accessible via 
internet \cite{AKcy}.
Work on toric complete intersections and nef partitions \cite{bb1,bb2,strh} 
is in progress. 
Further web pages with relevant information are \cite{RScy,SKcy}.

In the next section we give a self-contained exposition of our classification
algorithm and of the results in 3 and 4 dimensions. In section 3 we 
discuss the implications of these results for the geometry of 
toric K3 and Calabi--Yau hypersurfaces. In section 4 we explain the toric
realization of fibrations where both the fibered space and the fiber have
vanishing first Chern classes.
We discuss how weight systems can 
be used to encode such fibrations and how this is related to fibrations in
weighted projective spaces. 
We also provide an appendix with several tables that summarise some of our 
results.

\section{Classification of Reflexive Polyhedra}

In this section we 
give a self-contained exposition of our methods and
results on the classification of reflexive polyhedra, without reference to
toric geometry. 
Nevertheless, as we will see in the next section, some of the concepts used
here, in particular the concept of weights, have interpretations in terms
of geometry.

A polytope in $\IR^n$ is the convex hull of a finite set of points in 
$\IR^n$, and for our present purposes a polyhedron is the same thing as a 
polytope
(in particular, it is always bounded, which need not be true if a polyhedron
is defined as the intersection of a finite number of half spaces).

We will be interested in the case where we have a pair of lattices 
$M\simeq \IZ^n$ and $N=\mao{Hom}(M,\IZ)\simeq\IZ^n$ and their real extensions
$M_\IR\simeq \IR^n$ and $N_\IR\simeq \IR^n$.
A polyhedron $\D\subset M_\IR$ is called a lattice (or integer) polyhedron if
the vertices of $\D$ lie in $M$.

\Def A polytope $\D\subset \IR^n$ has the `interior point property' or `IP 
property', if $\ipo$ (the origin of $\IR^n$) is in the interior. 
A simplex with this property is an IP simplex.

\Def For any set $\D\subset M_\IR$ the dual 
(or polar) set $\D^*\subset N_\IR=M_\IR^*$ is given by
\beq 
\D^*=\{y\in N_{\IR}: ~~~\<y,x\>\ge -1 ~~~\forall x\in \D\},  
\eeql{dual}
where $\<y,x\>$ is the duality pairing between $y\in N_{\IR}$ and $x\in M_\IR$.

If $\D$ is a polytope with the IP property, then $\D^*$ is also a polytope 
with the IP property and $(\D^*)^*=\D$.

\Def A lattice polyhedron $\D\subset M_{\IR}$ is called reflexive if its dual 
$\D^*\subset N_{\IR}$ is a lattice polyhedron w.r.t. the lattice $N$ dual
to $M$.

The main idea of our classification scheme is to construct a set of polyhedra
such that every reflexive polyhedron is a subpolyhedron of one of the
polyhedra in this set. 
By duality, every reflexive polyhedron must contain one of the duals of these
polyhedra, so we are looking for polyhedra that are minimal in some sense.
In the following subsection we will give a definition of minimality that 
depends only on the way in which a polytope is spanned by its vertices, 
without reference to a lattice or details of the linear structure.
We will see that this allows for a very rough classification with only a few
objects in low dimensions.
The corresponding characterisation of polyhedra can be refined by specifying
explicitly the linear relations between the vertices with the help of weight 
systems.
We will see that these weight systems can be used in a simple way to find the
polyhedra dual to the minimal ones and to check whether they can possibly 
contain reflexive polyhedra;
the main criterion here is the existence of a dual pair of lattices such that
a minimal polytope is a lattice polyhedron and the convex hull of the 
lattice points of the dual has the IP property.
The classification of the relevant weight systems leads to a finite number of 
polytopes that contain all reflexive polytopes, with the subtlety that only
the linear structure but not the lattice on which some polytope may be 
reflexive is specified.
In the final subsection we solve this problem by showing how to 
identify all lattices on which a polyhedron given in terms of its linear 
structure can be reflexive, and present the results of our classification 
scheme.

\subsection{Minimal polyhedra and their structures}

We will later give various definitions of minimality, each of which has
advantages and disadvantages. 
Here we define the weakest form of minimality, but the one that is most useful,
where we forget for the time being about the lattice structure and concentrate
on the vertex structure only.

\Def A minimal polyhedron $\nabla\subset {\IR^n}$ is 
defined by the following properties:\\
1. $\nabla$ has the IP property.\\
2. If we remove one of the vertices of $\nabla$, the convex hull of the
  remaining vertices of $\nabla$ does not have the IP property.

Obviously every polytope $\nabla\subset {\IR^n}$ with the IP property contains
at least one minimal polytope spanned by a subset of the vertices of 
$\nabla$.
Before asking ourselves which minimal polytopes can be subpolytopes of 
reflexive polyhedra, we will now analyse the possible general structures of 
minimal polytopes.

\ni {\bf Lemma 1:} 
A minimal polytope $\nabla\subset\IR^n$ with vertices $V_1,\cdots,V_{k}$
is either a simplex or contains an $n'$-dimensional minimal
polytope $\nabla':=\mao{Convex Hull}\{V_1,\cdots,V_{k'}\}$ and an IP simplex
$S:=\mao{Convex Hull}(R\cup\{V_{k'+1},\cdots,V_{k}\})$ with
$R\subset\{V_1,\cdots,V_{k'}\}$ such that $k-k'=n-n'+1\ge 2$ and
$\mao{dim}S\le n'$.\\
{\it Proof:}
If $\nabla$ is a simplex, there is nothing left to prove. Otherwise,
we first note that every vertex $V$ of $\nabla$ must belong to at least
one IP simplex: It is always possible to find a triangulation of $\nabla$
such that every $n$-simplex in this triangulation has $V$ as a vertex
(just triangulate the cone whose apex is $V$ and whose one dimensional rays
are $V\tilde V$, where the $\tilde V$ are the other vertices of $\nabla$).
As $\ipo$ must belong to at least one of these simplices, it must lie on
some simplicial face which then is an IP simplex. Now
consider the set of all IP simplices consisting of vertices of $\nabla$. 
Any subset of this set will define a lower dimensional minimal polytope:
The fact that $\ipo$ is interior to each simplex means that it is a positive
linear combination of the vertices of any such simplex, and therefore $\ipo$ 
can also be written as a positive linear combination of all vertices involved.
If the corresponding polytope were not minimal, our original $\nabla$ could not
be minimal, either.
Among all lower dimensional minimal polytopes, 
take one (call it $\nabla'$) with the maximal dimension $n'$
smaller than $n$.
$\IR^n$ factorizes into $\IR^{n'}$ and $\IR^n/\IR^{n'}\cong\IR^{n-n'}$
(equivalence classes in $\IR^n$). 
The remaining vertices define a polytope $\nabla_{n-n'}$ in $\IR^n/\IR^{n'}$. 
If $\nabla_{n-n'}$ were not a simplex,
it would contain a simplex of dimension smaller than $n-n'$ which would
define, together with the vertices of $\nabla'$, a minimal polytope
of dimension $s$ with $n'<s<n$, in contradiction with our assumption.
Therefore $\nabla_{n-n'}$ is a simplex. 
Because of minimality of $\nabla$, each of the $n-n'+1$ vertices of 
$\nabla_{n-n'}$ can have only one representative
in $\IR^n$, implying $k-k'=n-n'+1$.
The equivalence class of $\ipo$ can be described uniquely as a positive
linear combination of these vertices.
This linear combination
defines a vector in $\IR^{n'}$, which can be written as a negative linear
combination of $\le n'$ linearly independent vertices of $\nabla'$.
These vertices, together with those of $\nabla_{n-n'}$, form the simplex $S$.
By the maximality assumption about $\nabla'$, dim$S$ cannot exceed 
dim$\nabla'$.
\hfill$\Box$

\Def For an n-dimensional minimal polytope $\nabla$ with $k$ vertices, an IP 
simplex structure is a collection of subsets $S_i,1\le i\le k-n$ of the set of 
vertices of $\nabla$, such that:\\
The convex hull of the vertices in each $S_i$ is an IP simplex, \\
$\nabla_j=\mao{Convex Hull}\bigcup_{i=1}^jS_i$ is a lower
dimensional minimal polytope for every $j\in\{1,\ldots,k-n\}$,\\
$\nabla_{k-n}=\nabla$ and\\
$S_j\setminus \bigcup_{i=1}^{j-1}S_i$ contains at least two vertices.

\ni {\bf Corollary:}
Every minimal polytope allows an IP simplex structure.\\
\Pro If $\nabla$ is a simplex, this is obvious. 
Otherwise one can choose $S_{k-n}=S$ and $\nabla_{k-n-1}=\nabla'$ with
$S$ and $\nabla'$ as in lemma 1 and proceed inductively. \hBo

\ni {\bf Lemma 2:} Denote by $\{S_i\}$ an IP simplex structure.
Then $S_i-\bigcup_{j\ne i}S_j$ never contains exactly one point.\\[2pt]
{\it Proof:} An IP simplex contains line
segments $VV'$ with $V'=-\e V$, where $\e$ is a positive number.
If a simplex $S=\mao{Convex Hull}\{V_1,\cdots, V_{s+1}\}$ has all of its
vertices except one ($V_{s+1}$) in common with
other simplices, then all points in the linear span of
$S$ are nonnegative linear combinations of the $V_j$ and the $-\e_jV_j$
with $j\le s$, thus showing that $V_{s+1}$ violates the minimality of 
$\nabla$.
\hfill$\Box$

The following example shows that an IP simplex structure need not be 
unique:\\[2mm]
{\bf Example:} $n=5$, $\nabla=\mao{Convex Hull}\{V_1,\cdots,V_8\}$ with
\bea V_1=(1,1,0,0,0),\;V_2=(1,-1,0,0,0),\;
     V_3=(-1,0,1,0,0),\;V_4=(-1,0,-1,0,0),\nn\\
     V_5=(-1,0,0,1,0),\;V_6=(-1,0,0,-1,0),\;
     V_7=(1,0,0,0,1),\;V_8=(1,0,0,0,-1).    \eea
$\nabla$ contains the IP simplices 
$S_{1234}=V_1V_2V_3V_4$ (in the $x_1x_2x_3$--plane),
$S_{1256}$ (in the $x_1x_2x_4$--plane),
$S_{3478}$ (in the $x_1x_3x_5$--plane),
$S_{5678}$ (in the $x_1x_4x_5$--plane) and the 4-dimensional
minimal polytopes
$\nabla_{123456}$, $\nabla_{123478}$, $\nabla_{125678}$, $\nabla_{345678}$.
Any set of three of the four IP simplices defines an IP simplex structure.

\ni {\bf Lemma 3:} 
For dimensions $n=1,2,3,4$ of $\IR^n$ precisely the
following IP simplex structures of minimal polyhedra are possible:\\
\begin{tabular}{ll}
$n=1$: & $\{S_1=V_1V_2\}$;\\
$n=2$: & $\{S_1=V_1V_2V_3\}$,\\
        & $\{S_1=V_1V_2,\; S_2=V_1'V_2',\}$;\\
$n=3$: & $\{S_1=V_1V_2V_3V_4\}$,\\
& $\{S_1=V_1V_2V_3,\; S_2=V_1'V_2'\}$,\\
& $\{S_1=V_1V_2V_3,\; S_2=V_1V_2'V_3'\}$,\\ 
& $\{S_1=V_1V_2,\; S_2=V_1'V_2',\; S_3=V_1''V_2''\}$;\\
$n=4$: & As in the first column of table 1 in the appendix.
\end{tabular}\\[2mm]
{\it Proof:} Recursive application of lemma 1 and use of lemma 2 shows that
these are the only possible structures. Explicit realisations of these
structures will be presented later.
\hfill$\Box$

\subsection{Weight systems}

Any IP polytope, and therefore any reflexive polyhedron, must obviously contain
one of the minimal polyhedra encountered in the last subsection. 
The structures found there are rather coarse, so now we have to face the task
of suitably refining them in such a way that they become useful for our goal
of classifying reflexive polyhedra. 
In particular, we will find that the linear relations between the vertices of
minimal polyhedra can be encoded by sets of real numbers called weight
systems, and we will adress the question of which weight systems can occur
if a minimal polyhedron is a subpolyhedron of some reflexive polytope.

The fact that a simplex spanned by vertices $V_i$ contains the origin in its
interior is equivalent to the condition that there exist positive real numbers
(weights) $q_i$ such that $\sum q_i V_i=0$. 
As these numbers are unique up to a common factor, it is convenient to choose
some normalization such as $\sum q_i=1$.

\Def A weight system is a collection of positive real numbers (weights) $q_i$
with $\sum q_i=1$.
A weight system corresponding to an IP simplex with vertices $V_i$ is the 
normalized set of numbers $q_i$ such that $\sum q_i V_i=0$.
A combined weight system (CWS) corresponding to a minimal polyhedron endowed 
with an IP simplex structure is the collection of weight systems $q_i^{(j)}$
corresponding to the IP simplices $S_j$ occurring there, with  $q_i^{(j)}=0$
if $V_i\not\in S_j$. 
We call a (combined) weight system rational if all of the $q_i$ are rational
numbers.

If a minimal polyhedron $\nabla$ is a lattice polyhedron, a
corresponding CWS will always be rational.
In this case it is possible to normalise the weights as positive integers
$n_i$ with no common divisor; then $q_i=n_i/d$ with $d=\sum n_i$.
We will use both conventions for describing weight systems.
By the definition of a lattice polyhedron, any lattice on which a minimal 
polyhedron $\nabla$ is integer must contain the lattice $N_{\rm coarsest}$ 
generated by the vertices of $\nabla$.

\Def Given a minimal polyhedron $\nabla\subset N_\IR$, we define the lattice 
$N_{\rm coarsest}$ as the lattice in $N_\IR$ generated linearly over $\IZ$ by 
the vertices of $\nabla$ and the lattice 
$M_{\rm finest}\subset M_\IR$ as the lattice dual to $N_{\rm coarsest}$.

\ni {\bf Lemma 4:} 
If $\nabla$ is a minimal polyhedron with vertices $V_i$ and $\q$ a CWS 
corresponding to an IP simplex structure of $\nabla$, then:\\
\del
a) With $\IR^k=\{(x_1,\ldots,x_k)\}$, the $k-n$ equations 
$\sum_i q_i^{(j)} x_i=0$ determined by
the weights corresponding to the IP simplices in $\nabla$ define a subspace
of $\IR^k$ isomorphic to $\IR^n$.\\
b) $\nabla^*$ is isomorphic to the polyhedron defined in this subspace by
$x_i\ge -1$ for $i=1,\ldots k$.\\
\enddel
a) The map $M_\IR\to \IR^k$, $X\to{\mathbf x}=(x_1,\ldots,x_k)$ with
$x_i=\<V_i,X\>$ defines an embedding such that the image of $M_\IR$ is 
the subspace defined by $\sum_i q_i^{(j)} x_i=0$ $\forall j$.\\
b) $\nabla^*$ is isomorphic to the polyhedron defined in this subspace by
$x_i\ge -1$ for $i=1,\ldots,k$.\\
c) If $\q$ is rational, then
$M_{\rm finest}$ is isomorphic to the sublattice of 
$\IZ^k=\{(x_1,\ldots,x_k) \hbox{ integer}\}\subset \IR^k$ determined by the 
equations $\sum_i q_i^{(j)} x_i=0$.\\[2pt]
\Pro 
a) $\sum q_i^{(j)} V_i=0$ implies $\sum_i q_i^{(j)} x_i=0$. 
Conversely, the $x_i$ determine $X$ because a
point in $M_\IR$ is uniquely determined by its duality pairings with a set
of generators (here, the $V_i$) of the dual space.\\
b) follows from the form of the embedding map and the definition of the
dual polytope (\ref{dual}).\\
c) If $X$ belongs to any lattice $M$ such that $\nabla$ is integer on the
dual lattice $N$, the corresponding $x_i$ must be integer. If the $x_i$ are
integer, then $X$ has integer pairings with the generators $V_i$ of 
$N_{\rm coarsest}$, so $X$ belongs to $M_{\rm finest}$.
\hBo

\ni {\bf Corollary:} 
An IP simplex structure together with the specification of a CWS uniquely 
determines a minimal polyhedron up to isomorphism.\\[2pt]
\Pro By lemma 4, $\nabla^*$ and hence $\nabla$ is uniquely
determined by the CWS.\hBo

As our example after lemma 2 shows, an IP simplex structure need not be
unique, so it is possible that two different CWS may correspond to the
same minimal polytope.
In such a situation, the weight systems of one CWS must be linear combinations
of those of the other CWS with coefficients that are not all nonnegative. 
Since all weights must be positive, this can only happen if there is an IP
simplex such that all of its vertices also belong to other IP simplices in
the same IP simplex structure.
This can happen only for $n\ge 5$, as one can see by explicitly checking
all cases for $n\le 4$. 
Thus, for $n\ge 5$ it might be preferable to work with equivalence classes 
of CWS leading to the same minimal polytopes instead of using CWS only.

\Def If $\q$ is a rational CWS corresponding to a minimal
polyhedron $\nabla$, we define $\D(\q)$ as the 
convex hull of $\nabla^*\cap M_{\rm finest}$.
We say that $\q$ has the IP property if $\D(\q)$ has the IP property.

\ni {\bf Corollary:} 
If a CWS has the IP property, then every single weight
system occurring in it also has the IP property.\\[2mm]
\Pro Without loss of generality we can assume that the single weight system is 
$\q^{(1)}$ with $q^{(1)}_i>0$ for $i\le l$ and $q^{(1)}_i=0$ for $i>l$.
There is a natural projection $\pi$ from $\IZ^k$ as in lemma 4 to $\IZ^l$
by restriction to the first $l$ coordinates.
Our construction implies that the projection of the lattice polytope
in $\IZ^k$ is a subpolytope of the lattice polytope in $\IZ^l$ determined 
by $\sum_i q_i^{(1)} x_i=0$ and $x_i\ge -1$.
If $\ipo^l=\pi(\ipo^k)$ were not in the interior of the polytope in $\IZ^l$, 
then $\ipo^k$ could not be in the interior of the polytope in $\IZ^k$.
\hBo

\del
\noindent
\begin{tabular}{||rrrr|r|c|} \hline\hline
$n_1$  & $n_2$  & $n_3$  & $n_4$  & d &   \\ \hline
1 & 1 & 1 & 1 & 4 & r\\[-6pt]
1 & 1 & 1 & 2 & 5 & r\\[-6pt]
1 & 1 & 2 & 2 & 6 & r\\[-6pt]
1 & 1 & 1 & 3 & 6 & r\\[-6pt]
1 & 1 & 2 & 3 & 7 & r\\[-6pt]
1 & 2 & 2 & 3 & 8 & s\\[-6pt]
1 & 1 & 2 & 4 & 8 & r\\[-6pt]
1 & 2 & 3 & 3 & 9 & s\\[-6pt]
1 & 1 & 3 & 4 & 9 & r\\[-6pt]
1 & 2 & 3 & 4 & 10 & s\\[-6pt]
1 & 2 & 2 & 5 & 10 & s\\[-6pt]
1 & 1 & 3 & 5 & 10 & r\\[-6pt]
1 & 2 & 3 & 5 & 11 & l\\[-6pt]
2 & 3 & 3 & 4 & 12 & s\\[-6pt]
1 & 3 & 4 & 4 & 12 & s\\[-6pt]
2 & 2 & 3 & 5 & 12 & s\\[-6pt]
1 & 2 & 4 & 5 & 12 & s\\[-6pt]
1 & 2 & 3 & 6 & 12 & s\\[-6pt]
1 & 1 & 4 & 6 & 12 & r\\[-6pt]
1 & 3 & 4 & 5 & 13 & l\\[-6pt]
2 & 3 & 4 & 5 & 14 & s\\[-6pt]
2 & 2 & 3 & 7 & 14 & s\\[-6pt]
1 & 2 & 4 & 7 & 14 & s\\[-6pt]
3 & 3 & 4 & 5 & 15 & s\\[-6pt]
2 & 3 & 5 & 5 & 15 & s\\[-6pt]
1 & 3 & 5 & 6 & 15 & s\\[-6pt]
1 & 3 & 4 & 7 & 15 & s\\[-6pt]
1 & 2 & 5 & 7 & 15 & s\\[-6pt]
1 & 4 & 5 & 6 & 16 & s\\[-6pt]
2 & 3 & 4 & 7 & 16 & s\\[-6pt]
1 & 3 & 4 & 8 & 16 & s\\[-6pt]
1 & 2 & 5 & 8 & 16 & s\\
\hline\hline \end{tabular}\hfil
\hbox{\begin{tabular}{|rrrr|r|c|} \hline\hline
$n_1$  & $n_2$  & $n_3$  & $n_4$  & d &   \\ \hline
2 & 3 & 5 & 7 & 17 & l\\[-6pt]
3 & 4 & 5 & 6 & 18 & s\\[-6pt]
1 & 4 & 6 & 7 & 18 & s\\[-6pt]
2 & 3 & 5 & 8 & 18 & s\\[-6pt]
2 & 3 & 4 & 9 & 18 & s\\[-6pt]
1 & 3 & 5 & 9 & 18 & s\\[-6pt]
1 & 2 & 6 & 9 & 18 & s\\[-6pt]
3 & 4 & 5 & 7 & 19 & l\\[-6pt]
2 & 5 & 6 & 7 & 20 & -\\[-6pt]
3 & 4 & 5 & 8 & 20 & s\\[-6pt]
2 & 4 & 5 & 9 & 20 & -\\[-6pt]
2 & 3 & 5 & 10 & 20 & s\\[-6pt]
1 & 4 & 5 & 10 & 20 & s\\[-6pt]
3 & 5 & 6 & 7 & 21 & -\\[-6pt]
1 & 5 & 7 & 8 & 21 & -\\[-6pt]
2 & 3 & 7 & 9 & 21 & s\\[-6pt]
1 & 3 & 7 & 10 & 21 & -\\[-6pt]
2 & 4 & 5 & 11 & 22 & -\\[-6pt]
1 & 4 & 6 & 11 & 22 & s\\[-6pt]
1 & 3 & 7 & 11 & 22 & -\\[-6pt]
3 & 6 & 7 & 8 & 24 & -\\[-6pt]
4 & 5 & 6 & 9 & 24 & s\\[-6pt]
1 & 6 & 8 & 9 & 24 & s\\[-6pt]
3 & 4 & 7 & 10 & 24 & s\\[-6pt]
2 & 3 & 8 & 11 & 24 & -\\[-6pt]
3 & 4 & 5 & 12 & 24 & s\\[-6pt]
2 & 3 & 7 & 12 & 24 & s\\[-6pt]
1 & 3 & 8 & 12 & 24 & s\\[-6pt]
4 & 5 & 7 & 9 & 25 & -\\[-6pt]
2 & 5 & 6 & 13 & 26 & -\\[-6pt]
1 & 5 & 7 & 13 & 26 & -\\[-6pt]
2 & 3 & 8 & 13 & 26 & -\\
\hline\hline \end{tabular}}\hfil
\hbox{\begin{tabular}{|rrrr|r|c||} \hline\hline
$n_1$  & $n_2$  & $n_3$  & $n_4$  & d &   \\ \hline
5 & 6 & 7 & 9 & 27 & -\\[-6pt]
2 & 5 & 9 & 11 & 27 & -\\[-6pt]
4 & 6 & 7 & 11 & 28 & -\\[-6pt]
3 & 4 & 7 & 14 & 28 & s\\[-6pt]
1 & 4 & 9 & 14 & 28 & -\\[-6pt]
5 & 6 & 8 & 11 & 30 & -\\[-6pt]
3 & 4 & 10 & 13 & 30 & -\\[-6pt]
4 & 5 & 6 & 15 & 30 & s\\[-6pt]
2 & 6 & 7 & 15 & 30 & -\\[-6pt]
1 & 6 & 8 & 15 & 30 & s\\[-6pt]
2 & 3 & 10 & 15 & 30 & s\\[-6pt]
1 & 4 & 10 & 15 & 30 & s\\[-6pt]
4 & 5 & 7 & 16 & 32 & -\\[-6pt]
2 & 5 & 9 & 16 & 32 & -\\[-6pt]
3 & 5 & 11 & 14 & 33 & -\\[-6pt]
4 & 6 & 7 & 17 & 34 & -\\[-6pt]
3 & 4 & 10 & 17 & 34 & -\\[-6pt]
7 & 8 & 9 & 12 & 36 & -\\[-6pt]
3 & 4 & 11 & 18 & 36 & -\\[-6pt]
1 & 5 & 12 & 18 & 36 & -\\[-6pt]
5 & 6 & 8 & 19 & 38 & -\\[-6pt]
3 & 5 & 11 & 19 & 38 & -\\[-6pt]
5 & 7 & 8 & 20 & 40 & -\\[-6pt]
3 & 4 & 14 & 21 & 42 & s\\[-6pt]
2 & 5 & 14 & 21 & 42 & -\\[-6pt]
1 & 6 & 14 & 21 & 42 & s\\[-6pt]
4 & 5 & 13 & 22 & 44 & -\\[-6pt]
3 & 5 & 16 & 24 & 48 & -\\[-6pt]
7 & 8 & 10 & 25 & 50 & -\\[-6pt]
4 & 5 & 18 & 27 & 54 & -\\[-6pt]
5 & 6 & 22 & 33 & 66 & -\\[-6pt]
 &  &  &  &  & \\
\hline\hline \end{tabular}}\hfill\\[3mm]
Table 2: IP weight systems for $l=4$. The last column indicates the minimality
type: `s' for span, `l' for lp-minimality and `r' for r-minimality.
As r-minimality implies lp-minimality and the latter implies the span property
for $l=4$, we have given only the strongest statement in each case.
\enddel

\ni {\bf Lemma 5:}
Let $l$ denote the number of weights of a weight system. Then the following
statements hold:\\[1mm]
$l=2$: There is a single IP weight system, namely $(1,1)$.\\
$l=3$: There are three IP weight systems, namely $(1,1,1)$, $(1,1,2)$ and 
         $(1,2,3)$. \\
$l=4$: There are the 95 IP weight systems shown in table 3.\\
$l=5$: There are 184,026 IP weight systems which
can be found at our web site \cite{KScy}.
\\[2mm]
\Pro
The classification of IP weight systems is based on the study of which integer
points are allowed by lemma 4.
Assume that a weight system $q_1,\cdots,q_l$ allows a collection of points with
coordinates $x_i\ge -1$ as in lemma 4, including the interior point with 
$x_i=0\;\forall i$.
If these points fulfill an equation of the type $\sum_{i=1}^la_ix_i=0$
with {\bf a}$\ne${\bf q}, then the weight system must also allow at least one
point with $\sum_{i=1}^la_ix_i>0$ and at least one point with
$\sum_{i=1}^la_ix_i<0$ to ensure that {\bf 0} is really in the interior.
The latter inequality is the one that we actually use for the algorithm:
Starting with the point {\bf 0}, we see that unless our
weight system is $\q=(1/l,\cdots,1/l)$, there must be at least one point
with $\sum_{i=1}^lx_i<0$. For $l\le5$ there are only a few possibilities,
and after choosing some point {\bf x}${}_1$, we can look for some simple
equation fulfilled by {\bf 0} and {\bf x}${}_1$ and proceed in the same way.
\\[2mm]
If $l=2$, any weight system except $(1/2,1/2)$ would have to allow an
integer point with $x_1+x_2<0$, $x_1\ge -1$ and $x_2\ge -1$. 
Such a point has no positive coordinate and therefore cannot be allowed by a
(positive) weight system.\\[2mm]
For $l=3$ the classification is still easily carried out by hand:
Unless $\q=(1/3,1/3,1/3)$, we need at least one point with $x_1+x_2+x_3<0$.
As points where no coordinate is greater than 0 would be in conflict with
the positivity of the weight system,
we need the point $(1,-1,-1)$ (up to a permutation of indices).
Now we note that {\bf 0} and $(1,-1,-1)$ both fulfill $2x_1+x_2+x_3=0$,
so $\q=(1/2,1/4,1/4)$ or we need a point with $2x_1+x_2+x_3<0$.
The only point allowed by this inequality which leads to a sensible
weight system is $(-1,2,-1)$, leading to $\q=(1/2,1/3,1/6)$.
\\[2mm]
For $l=4$ and $l=5$ we have implemented this strategy in a computer program
that produced 99 and 200653 candidates for IP weight systems, respectively.
Finally, explicit constructions of $\D(\q)$ show that four of the 99 weight
systems with $l=4$ and 16627 of the 200653 weight systems with $l=5$ do not
have the IP property, leading to the results given.
\hfill$\Box$

\ni {\bf Remark:} The 95 IP weight systems for $l=4$ are precisely the well
known 95 weight systems for weighted $\IP^4$'s that have K3 hypersurfaces 
\cite{reid,fl89}, whereas for $l=5$ the 7555 weight systems corresponding to 
weighted $\IP^4$'s that allow transverse polynomials \cite{KlS,nms} are just 
a small subset of the 184026 different IP weight systems.

\ni {\bf Lemma 6:} 
In dimensions $n=1,2,3,4$, the CWS with the IP property
are the weight systems with $l=n+1$ given in the previous lemma and, in 
addition, the following CWS:\\
$n=2$: $\{(1,1,0,0),(0,0,1,1)\}$ \\
$n=3$: The 21 CWS given in table 2\\
$n=4$: 17320 CWS (cf. the second column of table 1)\\[2mm]
\Pro By explicitly combining the structures of Lemma 3 with the IP weight
systems of Lemma 5 and checking for the IP property of $\D(\q)$. \hBo

\subsection{The classification}

As we saw in the previous subsections, every reflexive polyhedron must contain
at least one minimal polytope corresponding to one of the CWS found there.
Thus, by duality, every reflexive polyhedron must be a subpolyhedron of one
of the $\D(\q)$ on some suitable sublattice of the finest possible lattice
$M_{\rm finest}$. 
We start this section with analysing the question of which dual pairs of
lattices can be chosen such that a dual pair of polyhedra is reflexive on
them.
Then we give various refinements of our original definition of minimality, and
finally we present our results on the classification of reflexive polyhedra.

Given a dual pair of polytopes such that $\D$ has $n_V$ vertices and $n_F$
facets (a facet being a codimension 1 face), the dual polytope has 
$n_V$ facets and $n_F$ vertices.

\Def The vertex pairing matrix (VPM) $X$ is the $n_F\times n_V$
matrix whose entries are $X_{ij}=\<\bar V_i,V_j\>$, where $\bar V_i$ and
$V_j$ are the vertices of $\D^*$ and $\D$, respectively.

$X_{ij}$ will be $-1$ whenever $V_j$ lies on the $i$'th facet.
Note that $X$ is independent of the choice of a dual pair of bases in 
$N_{\IR}$ and $M_\IR$ but depends on the orderings of the vertices.
If $\D$ is reflexive, then its VPM is obviously integer.
In this case there are distinguished lattices $M_{\rm coarsest}$ and 
$N_{\rm coarsest}$, generated by the vertices of $\D$ and $\D^*$, respectively,
and their duals $N_{\rm finest}$ and $M_{\rm finest}$.
Clearly any lattice $M$ on which $\D$ is reflexive must fulfill
$M_{\rm coarsest}\subseteq M\subseteq M_{\rm finest}$.

\ni {\bf Lemma 7:} 
If $\D\subset M_\IR\simeq\IR^n$ is a polytope with the IP property such that
its VPM $X$ is integer, the following statements hold:\\
$X$ can be decomposed as 
$X=\tilde W\cdot \tilde D\cdot \tilde U= W\cdot D\cdot U$,
where $\tilde W$ is a $GL(n_F,\IZ)$ matrix, $\tilde U$
is a $GL(n_V,\IZ)$ matrix and $\tilde D$ is an $n_F\times n_V$
matrix such that the first $n$
diagonal elements are positive integers whereas all other elements are zero;
$W$, $D$ and $U$ are the obvious $n_F\times n$, $n\times n$ and $n\times n_V$
submatrices.\\
The lattices $M\subset M_\IR$ on which $\D$ is reflexive are in one to one 
correspondence with decompositions $D=T\cdot S$, where $T$ and 
$S$ are upper triangular integer matrices with positive diagonal elements
and with $0\le T_{ji}<T_{ii}$. 
Then $\D$ as a lattice polyhedron on $M$ is isomorphic to the polytope in
$\IZ^n$ whose vertices are given by the columns of $S\cdot U$ 
and $\D^*$ is isomorphic to the polytope in $\IZ^n$ whose vertices are given 
by the lines of $W\cdot T$.
In particular, $\D$ on $M_{\rm finest}$ corresponds to 
$D\cdot U$, $\D$ on $M_{\rm coarsest}$ corresponds 
to $U$, $\D^*$ on $N_{\rm finest}$ corresponds to 
$W\cdot D$ and $\D^*$ on $N_{\rm coarsest}$ 
corresponds to $W$.
\\[2pt]
\Pro By recombining the lines and columns of $X$ in the style of Gauss's
algorithm for solving systems of linear equations, we can
turn $X$ into an $n_F\times n_V$ matrix $\tilde D$ with non-vanishing 
elements only along the diagonal.
But recombining lines just corresponds to left multiplication with 
some $GL(n_F,\IZ)$ matrix, whereas recombining columns corresponds to right
multiplication with some $GL(n_V,\IZ)$ matrix.
Keeping track of the inverses of these matrices, we successively create
decompositions $X=\tilde W^{(n)}\cdot \tilde D^{(n)}\cdot \tilde U^{(n)}$ 
(with $\tilde W^{(0)}=1$, $\tilde D^{(0)}=X$ and $\tilde U^{(0)}=1$).
We denote the matrices resulting from the last step by 
$\tilde W$, $\tilde D$ and $\tilde U$.
$\tilde W$ and $\tilde U$ being regular matrices and the rank of $X$ being
$n$, it is clear that $\tilde D$ has only $n$ non-vanishing elements
which can be taken to be the first $n$ diagonal elements.

In the same way as we defined an embedding of $M_\IR$ in $\IR^k$ lemma 4,
we now define an embedding in $\IR^{n_F}$ such that
$M_{\rm finest}$ is isomorphic to the sublattice of $\IZ^{n_F}$ 
determined by the linear relations among the $\bar V_i$.
In this context the $X_{ij}$ are just the embedding coordinates of the $V_j$.
The $n_F\times n_F$ matrix $\tilde W$ effects a change of coordinates in 
$\IZ^{n_F}$ so that $\D$ now lies in
the lattice spanned by the first $d$ coordinates.
Thus we can interpret the columns of $D\cdot U$ as the vertices of $\D$ on 
$M_{\rm finest}$.
Similarly, the lines of $W\cdot D$ are coordinates of the vertices of 
$\D^*$ on $N_{\rm finest}$, whereas $U$ and $W$ are the corresponding 
coordinates on the coarsest possible lattices.

Denoting the generators of $M_{\rm coarsest}$ by $\1 E_i$ 
and the generators of $M_{\rm finest}$ by $\1e_i$, we have 
$\1 E_i=\1e_jD_{ji}$.
An intermediate lattice will have generators $\1\ce_i=\1e_jT_{ji}$
such that the $\1 E_i$ can be expressed in terms of the $\1\ce_j$,
amounting to $\1 E_i=\1\ce_jS_{ji}=\1e_kT_{kj}S_{ji}$
with some integer matrix $S$.
This results in the condition $D_{ki}=T_{kj}S_{ji}$.
In order to get rid of the redundancy coming from the fact that the 
intermediate lattices can be described by different sets of generators,
one may proceed in the following way: 
$\1\ce_1$ may be chosen as a multiple of $\1e_1$ (i.e., 
$\1\ce_1=\1e_1T_{11}$).
Then we choose $\1\ce_2$ as a vector in the $\1e_1$-$\1e_2$-plane 
(i.e., $\1\ce_1=\1e_1T_{12}+\1e_2T_{22}$) subject to the
condition that the lattice generated by $\1\ce_1$ and $\1\ce_2$ should be 
a sublattice of the one generated by $\1E_1$ and $\1E_2$, which is 
equivalent to the possibility of solving $T_{kj}S_{ji}=D_{ki}$ for integer 
matrix elements of $S$.
We may avoid the ambiguity arising by the possibility of adding a multiple of
$\ce_1$ to $\ce_2$ by demanding $0\le T_{12}<T_{11}$.
We can choose the elements of $T$ column by column (in rising order).
For each particular column $i$ we first pick $T_{ii}$ such that
it divides $D_{ii}$; then $S_{ii}=D_{ii}/T_{ii}$.
Then we pick the $T_{ji}$ with $j$ decreasing from $i-1$ to $1$.
At each step the $j$'th line of $T\cdot S =D$,
\beq T_{ji}S_{ii}+\sum_{j<k<i}T_{jk}S_{ki}+T_{jj}S_{ji}=0, \eeq
must be solved for the unknown $T_{ji}$ and $S_{ji}$ with the
extra condition $0\le T_{ji}<T_{ii}$ ensuring that we get only one
representative of each equivalence class of bases.
\hBo

At this point we have, in principle, all the ingredients that we need for a 
complete classification of reflexive polyhedra.
We simply have to construct all subpolyhedra with integer VPM of all $\D(\q)$ 
with $\q$ being one of our IP CWS, and apply lemma 7.
Both for theoretical and for practical reasons, however, it is interesting to 
reduce the number of polyhedra used as a starting point in our scheme. 
To this end we will give various refinements of our original definition of
minimality, preceded by a useful lemma on the structure of $\D(\q)$.

\ni {\bf Lemma 8:} 
For $n\le 4$, $\D(\q)$ is reflexive whenever it has the IP property.\\[2pt]
\Pro This fact was proved in \cite{wtc} and later explicitly confirmed by our
computer programs. \hBo

\Def Let $\nabla\subset N_{\IR}$ be a minimal lattice polyhedron 
such that $\D$, the convex hull of $\nabla^*\cap M$, also has the IP property.
Then we say that\\
$\nabla$ has the span property if the vertices of $\nabla$ are also vertices
of $\D^*$.\\
$\nabla$ is lp-minimal: If we remove one of the vertices of $\nabla$, the 
   convex hull of the remaining set of {\sl lattice points} of $\nabla$ does 
   not have the IP property.\\
$\nabla$ is very minimal: 
   If we remove one of the vertices of $\nabla$ from 
   the set of lattice points of $\D^*$, the convex hull of the remaining 
   lattice points of $\D^*$ does not have the IP property.\\
A CWS $\q$ is said to have one of the above properties if the corresponding
$\nabla$ on $N_{\rm coarsest}$ has it.\\
A reflexive polytope $\D\subset M_\IR$ is called r-maximal (and its dual 
$\D^*\subset N_\IR$  r-minimal) if it is not contained in any other reflexive 
polytope. \\
A CWS $\q$ is called r-minimal if $\D(\q)$ is r-maximal.

The name `span property' refers to the fact that our definition is equivalent 
to the statement that the hyperplanes in $M_\IR$ dual to the vertices of 
$\nabla$ are spanned by points of $\D$.
The following lemma clarifies the relations between the various definitions
of minimality and the ways in which these definitions can be used to refine 
our classification scheme. 
It also answers the question of how many CWS of the various minimality types 
exist.

\ni {\bf Lemma 9:} \\
a) For every reflexive polytope $\D\subset M_\IR$, there exists at least one
CWS $\q$ with the span property such that $\D$ is a subpolyhedron of 
the convex hull of $\nabla^*\cap M$ and $M$ is a sublattice of 
$M_{\rm finest}$.\\
b) For every reflexive polytope $\D\subset M_\IR$, there exists at least one
lp-minimal CWS $\q$ such that $\D$ is a subpolyhedron of the convex hull of 
$\nabla^*\cap M$ and $M$ is a sublattice of $M_{\rm finest}$.\\
c) If $\q$ is very minimal, $\D(\q)$ is not a subpolyhedron of $\D(\q')$ for 
any $\q'$ corresponding to a minimal polytope different from the one defined 
by $\q$.\\
d) A very minimal polytope is lp-minimal and has the span property.\\
e) For every reflexive polytope $\D\subset M_\IR\simeq\IR^n$ with $n\le 4$, 
there exists at least one r-minimal CWS $\q$ such that $\D$ is a subpolyhedron 
of the convex hull of $\nabla^*\cap M$ and $M$ is a sublattice of 
$M_{\rm finest}$.\\
f) For $n\le 4$, a CWS $\q$ is r-minimal if and only if it is very minimal.\\
g) For $n\le 3$ (but not for $n=4$), every lp-minimal CWS has the span 
property.\\
h) The very minimal CWS for $n=2$ are $\{(1,1,1)\}, \{(1,1,2)\}$ and 
$\{(1,1,0,0),(0,0,1,1)\}$. The remaining IP weight system $\{(1,2,3)\}$ has
the span property but is not lp-minimal.\\
i) For $n=3$ the minimality type is indicated in tables 2 and 3. \\
j) For $n=4$ the numbers of CWS of the different minimality types are given in
table 1.
\\[2mm]
\Pro a) By dropping vertices from $\D^*$ one can always arrive at a minimal
polytope $\nabla\subseteq\D^*$ and the corresponding CWS.\\
b) By dropping lattice points from $\D^*$ one can always arrive at an 
lp-minimal (and therefore also minimal) polytope $\nabla\subseteq\D^*$ and 
the corresponding CWS.\\
c) If $\D(\q)$ were a subpolyhedron of $\D(\q')$ for some $\q'$ other than 
$\q$, then $(\D(\q))^*$ would contain (but not be equal to) $(\D(\q'))^*$, 
which is impossible by the definition of $\q$ being very minimal.\\
d) By definition, $\D\subseteq \nabla^*$, implying $\nabla \subseteq \D^*$.
Very minimal implies span: If a vertex of $\nabla$ were not a vertex 
of $\D^*$, it would be in the convex hull of the remaining lattice points 
of $\D^*$ which then would be equal to $\D^*$ and hence have the IP property,
thus violating the assumption that $\nabla$ is very minimal.
The fact that very minimal implies lp-minimal is obvious from comparing the 
different definitions.\\
e) With a), we can find a CWS $\q^{(1)}$ such that $\D$ is a subpolyhedron 
of $\D(\q^{(1)})$, possibly on a sublattice. By lemma 8, $\D(\q^{(1)})$ is
reflexive. If $\q^{(1)}$ is not r-minimal, $\D(\q^{(1)})$ is a proper 
subpolyhedron of
some other reflexive polyhedron $\D^{(1)}$ for which we can find a CWS 
$\q^{(2)}$ as before. As the number of lattice points of $\D(\q^{(i)})$ 
increases in every step, this process has to terminate; thus $\q^{(i)}$ 
must be r-minimal for some $i$.\\
f) Because of c), every very minimal CWS is r-minimal. The fact that every 
r-minimal CWS is very minimal was checked explicitly by our computer 
programs.\\
g) -- j) By explicit checks, for $n\ge 3$ with the help of our computer 
programs.
\hBo

To end this section, we now give the results of the application of our 
classification scheme for various dimensions.

\ni {\bf Proposition 1:} 
For $n=2$ there are 16 reflexive polyhedra up to linear isomorphisms.
All of them are subpolyhedra of $\D(\q)$ where $\q$ is one of the three very 
minimal CWS.\\[2mm]
\Pro The classification of 2-dimensional reflexive polyhedra has been
established for a while (see, e.g., \cite{ba85,ko90}) and is easily reproduced
within our scheme. The second fact can be checked explicitly. \hBo

\ni {\bf Proposition 2:} 
For $n=3$ there are 4319 reflexive polyhedra up to linear isomorphisms.
4318 of them are subpolyhedra of $\D(\q)$ where $\q$ is one of the very 
minimal CWS of tables 2 and 3. 
The remaining one is the convex hull of 
$\nabla^*\cap M$, where $\nabla$ is determined by the weight system
$(1,1,1,1)$ and $M$ is a $\IZ_2$ sublattice of $M_{\rm finest}$.\\[2mm]
\Pro We explicitly constructed all subpolyhedra with integer VPMs of the 
$\D(\q)$ coming from very minimal CWS with the help of a computer program and 
checked that polyhedra coming from CWS that are not very
minimal are contained in the list of reflexive subpolyhedra of the $\D(\q)$
for very minimal CWS.
Application of lemma 7 produced the last polyhedron. \hBo

\ni {\bf Proposition 3:} 
For $n=4$ there are more than 473.8 million  
reflexive polyhedra up to linear isomorphisms.
In addition to $\D(\q)$ with $\q$ one of 308 the r-minimal CWS, there are at 
least 25 further r-maximal polyhedra. \\[2mm]
\Pro Our computer programs have already produced 
more than 473.8 million different reflexive polyhedra. 
The 25 additional r-maximal polyhedra were obtained by applying lemma 7 to the
original 308 r-maximal polytopes and checking for r-minimality of the duals on
the various lattices allowed by lemma 7. \hBo

\section{Geometric interpretation of our classification results}

We now want to discuss what our results on the classification of reflexive 
polyhedra imply for Calabi--Yau manifolds that are hypersurfaces in toric 
varieties.

The lattice points of a reflexive polyhedron $\D$ encode the monomials
occurring in the description of the hypersurface in a variety $V_\S$ whose
fan $\S$ is determined by a triangulation of the dual polyhedron $\D^*$.
For details of what a fan is and how it determines a toric variety,
it is best to look up a standard textbook \cite{Ful,Oda}.
There is one particular approach to the description of toric varieties, however,
which cannot be found there.
This is the description in terms of  homogeneous coordinates \cite{Cox1},
which is the one most useful for applications in physics, and which also 
exhibits in the clearest way the significance of the weight systems that we
used in the context of our classification scheme.
We will briefly present this approach and show how Calabi--Yau manifolds are
constructed in this setup and then we will proceed to explain some of the 
consequences of our results in terms of geometry.

Given a fan $\S$ in $N_\IR$, it is possible to assign
a global homogeneous coordinate system to $V_\S$ in a way
similar to the usual construction of $\IP^n$.
To this end one assigns a coordinate $z_k$, $k=1,\cdots,K$
to each one dimensional cone in $\S$. 
If the primitive generators $v_1,\cdots,v_K$ of these one 
dimensional cones span $N_\IR$,
then there must be $K-n$ independent linear relations of
the type $\sum_k w^k_j v_k=0$. 
These linear relations are used to define equivalence relations
of the type 
\beq (z_1,\cdots,z_K)\sim (\l^{w^1_j}z_1,\cdots,\l^{w^K_j}z_K),~~~~~~
j=1,\cdots,K-n  \eeql{er} 
on the space $\IC^K\setminus Z_\S$. 
The set $Z_\S$ is determined by the fan $\S$ in the following way:
It is the union of spaces 
$\{(z_1,\cdots,z_K):\; z_i=0\;\forall i\in I\}$,
where the index sets $I$ are those sets for which $\{v_i:\;i\in I\}$
does not belong to a cone in $\S$.
Thus $(\IC^*)^K\subset \IC^K \setminus Z_\S \subset \IC^K \setminus \{0\}$.
Then $V_\S=(\IC^K \setminus Z_\S)/((\IC^*)^{(K-n)}\times G)$, where the $K-n$ 
copies of $\IC^*$ act by the equivalence relations given above and the finite 
abelian group $G$ is the quotient of the $N$ lattice by the lattice generated 
by the $v_k$.
We will usually consider the case where $G$ is trivial.
In this approach the toric divisors $D_k$ are determined by the equations
$z_k=0$.

The construction of a Calabi--Yau
hypersurface from a reflexive polyhedron proceeds in the following way:
We take $\D$ to be a reflexive polyhedron in $M_\IR$,
$\D^*\subset N_\IR$ its dual, and $\S$ a fan defined by
a maximal triangulation of $\D^*$.
This means that the integer generators $v_1,\cdots,v_K$ of the one 
dimensional cones are just the integer points (except the origin)
of $\D^*$. 
The polynomial whose vanishing determines the Calabi--Yau hypersurface 
takes the form
\beq \sum_{x\in \D\cap M}a_x\prod_{k=1}^K z_k^{\<v_k,x\>+1}.  \eeql{polyn}
It is easily checked that it is quasihomogeneous with respect to all 
$K-n$ relations of (\ref{er}) with degrees $d_j=\sum_{k=1}^K w_j^k$,
$j=1,\cdots K-n$.
Note how the reflexivity of the polyhedron ensures that the
exponents are nonnegative.

By \cite{bat}, the Hodge numbers $h_{11}$ and $h_{1,n-2}$ are known,
and in \cite{BD} the remaining Hodge numbers of the type $h_{1i}$
were calculated. 
For a hypersurface of dimension $n-1\ge 2$ these formulas can be summarised as
\bea h_{1i}&=&\d_{1i}\(l(\D^*)-n-1-\sum_{{\rm codim }\th^* =1}l^*(\th^*)\)+
        \d_{n-2,i}\(l(\D)-n-1-\sum_{{\rm codim }\th =1}l^*(\th)\)\nn\\
           &&+ \sum_{{\rm codim }\th^* =i+1}l^*(\th^*)l^*(\th) \label{baho}
                                                                       \eea
for $1\le i\le n-2$, where $l$ denotes the number of integer points of a 
polyhedron and $l^*$ denotes the number of interior integer points of a face.
These formulas are invariant under the simultaneous
exchange of $\D$ with $\D^*$ and $h_{1i}$ with $h_{1,n-i}$
so that Batyrev's construction is manifestly mirror symmetric
(at least at the level of Hodge numbers).
For $n\le 4$, the generic $(n-1)$-dimensional Calabi-Yau hypersurface in the
family defined by $\D$ will be smooth \cite{bat} and the meaning of these 
numbers is unambiguous.
For $n\ge 5$, the Calabi-Yau variety may have singularities that do not 
allow a crepant blow-up.
In this case we refer the reader to refs. \cite{BD} for a discussion of
the precise meaning of the Hodge numbers resulting from eq. (\ref{baho}).

In the case of a K3 surface there is only one such number, namely $h_{11}$,
which is well known always to be equal to 20. 
Contrary to the case of higher dimensional Calabi-Yau manifolds,
this number is not the same as the Picard number, which is given by \cite{bat}
\beq \hbox{Pic}=l(\D^*)-4-\sum_{{\rm facets}~\th^*~{\rm of}~\D^*}l^*(\th^*)+
     \sum_{{\rm edges}~\th^*~{\rm of}~\D^* }l^*(\th^*)l^*(\th).
\eeq
Mirror symmetry for K3 surfaces is usually interpreted in terms of 
families of lattice polarized K3 surfaces 
(see, e.g., \cite{Dol} or \cite{Arev}).
In this context the Picard number of a generic element of a family and the 
Picard number of a generic element of the mirror family add up to 20.
The fact that the Picard numbers for toric mirror families add up to
$20+\sum l^*(\th^*)l^*(\th)$ indicates that our toric models
occupy rather special loci in the total moduli spaces.

If a polyhedron $\D_1$ contains a polyhedron $\D_2$, then the definition
of duality implies $\D^*_1\subset \D^*_2$. 
Therefore the variety determined by the fan over $\D^*_1$ may be obtained
from the variety determined by the fan over $\D^*_2$ by blowing down
one or several divisors.
If we perform this blow-down while keeping the same monomials (those 
determined by $\D_2$), we obtain a generically singular hypersurface.
This hypersurface can be desingularised by varying the complex structure 
in such a way that we now allow monomials determined by $\D_1$.
Thus the classes of Calabi-Yau hypersurfaces determined by polyhedra $\D_1$ and
$\D_2$, respectively, can be said to be connected whenever $\D_1$ contains
$\D_2$ or vice versa.
More generally, if there is a chain of polyhedra $\D_i$ such that $\D_i$
and $\D_{i+1}$ are connected in the sense defined above, we call the
hypersurfaces corresponding to any two elements of the chain connected.

We can easily check for connectedness as a by-product of our classification 
scheme: 
For each new CWS $\q$ we check explicitly that at least one of the 
subpolyhedra of $\D(\q)$ has been found before. 
Connectedness of the corresponding list of 4318 polytopes in three dimensions
follows from the fact that this is always the case.
Connectedness of all 3-d reflexive polyhedra follows from the fact that
the last polytope that we only found on a sublattice contains 679 reflexive 
proper subpolytopes that were found before.
In the same way all of the four dimensional polyhedra that we have found so far
form a connected web.

As we saw in the previous section, every three dimensional reflexive polytope
$\D^*\subset N_\IR$ contains one of 16 r-minimal polyhedra as a subpolytope
on the same lattice.
Therefore, the fan of any toric ambient variety determined by a maximal
triangulation of a reflexive polyhedron is a refinement of one of the 
corresponding 16 fans.
In other words, any such toric ambient variety is given by the blow-up 
of one of the following 16 spaces  (cf. tables 2 and 3 and proposition 2):\\
-- $\IP^3$, \\
-- $\IP^3/\IZ_2$, \\
-- 8 different weighted projective spaces $\IP^2_{(q_1,q_2,q_3)}$,\\
-- $\IP^2\times\IP^1$,  \\
-- $\IP^2_{(1,1,2)}\times\IP^1$, \\
-- 3 further double weighted spaces, and\\
-- $\IP^1\times\IP^1\times\IP^1$. \\
Each of the three spaces with `overlapping weights' allows two distinct
bundle structures: The first one can be interpreted as a $\IP^2$ bundle
in two distinct ways, the second one as a $\IP^2$ bundle or a 
$\IP^2_{(1,1,2)}$ bundle, and the third one can be interpreted as a 
$\IP^2_{(1,1,2)}$ bundle in two distinct ways.
In each case the base space is $\IP^1$.

Let us end this section with briefly discussing a few of the most interesting
objects in our lists.
There are precisely two mirror pairs with Picard numbers
1 and 19, respectively.
One of them is the quartic hypersurface in $\IP^3$ with Picard number 1,
together with its mirror of Picard number 19, which is also the model 
whose Newton polytope is the only reflexive polytope with only 5 
lattice points.
This model corresponds to a blow-up of a $\IZ^4\times\IZ^4$ orbifold
of $\IP^3$.
The blow-up of six fixed lines $z_i=z_j$ by three divisors each yields
18 exceptional divisors leading to the total Picard number of 19.
The other mirror pair with Picard numbers 1 and 19 consists of the hypersurface
in $\IP^3_{(1,1,1,3)}$ of degree 6 and an orbifold of the same model, with
Newton polyhedra with 39 and 6 points, respectively.
This polyhedron is also one of the two `largest' polyhedra in the sense that 
there is
no reflexive polytope in three dimensions with more than 39 lattice points.
The other polyhedron with the maximal number of 39 points is the Newton 
polytope of the hypersurface of degree 12 in $\IP^3_{(1,1,4,6)}$.
This model leads to the description of elliptically fibered K3 surfaces
that is commonly used in F-theory applications \cite{Fv,MVI,MVII},
with the elliptic fiber embedded in a $\IP^2_{(1,2,3)}$ by a Weierstrass
equation.
The mirror family of this class of models can be obtained by forcing two
$E_8$ singularities into the Weierstrass model and blowing them up.
The resulting hypersurface allows also a different fibration structure
which can develop an $SO(32)$ singularity; thereby this model is able
to describe the F-theory duals of both the $E_8\times E_8$ and the $SO(32)$
heterotic strings with unbroken gauge groups in 8 dimensions \cite{fst}.

In four dimensions there is a unique `largest' object, determined by
the weight system $(1,1,12,28,42)/84$.
It has the maximum number, namely 680, of lattice points and the corresponding
Calabi--Yau threefold has the Hodge numbers $h_{11}=11$ and $h_{12}=491$.
The latter is the largest single Hodge number in our list, 
and the value of $|\chi|=|2(h_{11}-h_{12})|=960$ is also maximal,
the only other object with the same values being the mirror.
F-theory compactifications of the latter lead to the theories with the
largest known gauge groups in six dimensions \cite{fst}.

Another interesting object that we encountered is the 24-cell,
a self dual polytope with 24 vertices, which leads to a self mirror 
Calabi--Yau manifold with Hodge numbers (20,20).
It has the maximal symmetry order $1152=24*48$ among all 4 dimensional
        reflexive polytopes and arises as a subpolytope of the hypercube.
        It is a Platonic solid that contains the Archimedian cuboctahedron 
        (with symmetry order 48) as a reflexive section through the origin
        parallel to one of its 24 bounding octahedra. Note that in our context
        symmetries are realized as lattice isomorphisms, i.e. as subgroups of
        $GL(n,\IZ)$, and not as rotations.

The polytope with the largest order, namely 128, of 
$M_{\rm finest}/M_{\rm coarsest}$ is determined by the weight system 
$(1,1,1,1,4)/8$. For the 
Newton polytope of the quintic hypersurface in $\IP^4$, this order is 125.
There is a well known $\IZ_5$ orbifold of the quintic 
with Hodge numbers (1,21) which is quite peculiar from the lattice point of 
view: Although the $N$ lattice is not the lattice $N_{\rm finest}$ generated 
by the vertices of $\D^*$, the only lattice points of $\D^*$ are its vertices
and the IP. 
Thus it provides an example where the $N$ lattice is not even generated
by the lattice points of $\D^*$. 
This can only happen in more than 3 dimensions:
As a lattice triangle with 3 lattice points is always regular (i.e. it has 
the minimal volume 1 in lattice units) and there are no lattice hyperplanes 
between a facet and the IP because of reflexivity, the vertices of any 
triangle of a maximal triangulation of a 2-dimensional facet of a 
3-dimensional polytope provide a lattice basis.

\section{Fibrations}

In this section we want to discuss
fibrations of hypersurfaces of holonomy $SU(n-1)$ in $n$-dimensional
toric varieties where the generic fiber is an ($n\fib-1$)--dimensional
variety of holonomy $SU(n\fib-1)$.
In other words, it will apply to elliptic fibrations of
K3 surfaces, CY threefolds, CY fourfolds, etc., to K3 fibrations of
CY $k$-folds with $k\ge 3$, to threefold fibrations of fourfolds, and so on.
The main message is that the structures occurring in the fibration
are reflected in structures in the $N$ lattice: 
The fiber, being an algebraic subvariety of the whole space,
is encoded by a polyhedron $\D^*\fib$ which is a subpolyhedron of
$\D^*$, whereas the base, which is a projection of the 
fibration along the fiber, can be seen by projecting the $N$ lattice
along the linear space spanned by $\D^*\fib$.
We will first give a general discussion and then explain how descriptions
in terms of CWS may be useful for identifying and/or encoding fibration 
structures.

\subsection{Fibrations and reflexive polyhedra} 

Assume that $\D^*$ contains a lower-dimensional reflexive 
subpolyhedron $\D^*\fib=(N\fib)_\IR\cap \D^*$ with the same interior
point.
This allows us to define a dual pair of exact sequences
\beq 0\to N\fib \to N\to N\bas \to 0\eeql{esN}
and
\beq 0 \to M\bas \to M \to M\fib \to 0,\eeql{esM}
and corresponding sequences for the underlying real vector spaces.
We can convince ourselves that the image of $\D$ under 
$M_\IR \to (M\fib)_\IR$ is dual to $\D^*\fib$ in the following way:
We choose a basis $e^j$, $j=1,\ldots n$ of $N$ such that
$N\fib$ is generated by the $e^j$ with $1\le j\le n\fib$
and define $e_i$ to be the dual basis.
Then
\bea
\D\fib&=&\{(x^1,\cdots,x^{n\fib}):\exists~x^{n\fib+1},\ldots,x^{n}
 \hbox{ with }
      (x^1,\cdots,x^{n})\in\D)\},\\
(\D^*)\fib&=&\{(y_1,\cdots,y_{n\fib}): (y_1,\cdots,y_{n\fib},0)\in\D^*\},
\eea
and the duality of these two polytopes is easily checked.

Let us also assume that the image $\S\bas$ of
$\S$ under $\p: N\to N\bas$ defines a fan in $N\bas$.
This is certainly not true for  arbitrary triangulations of $\D^*$. 
Constructing fibrations, one should rather build a fan $\S\bas$
from the images of the one-dimensional cones in $\S$
and try to construct a triangulation of $\S$ and thereby of
$\D^*$ that is compatible with the projection.
It would be interesting to know whether this is always possible
whenever the intersection of a reflexive polyhedron with a linear 
subspace of $N_\IR$ is again reflexive.

The set of one-dimensional cones in $\S\bas$ is the set of 
images of one-dimensional cones in $\S$ that do not lie in $N\fib$. 
The image of a primitive generator $v_i$
of a cone in $\S$  is the origin or
a positive integer multiple of a primitive
generator $\tilde v_j$ of a one-dimensional cone in $\S\bas$. 
Thus we can define a matrix $r^i_j$, most of whose elements are $0$, through
$\p v_i=r_i^j\tilde v_j$ with $r_i^j\in \IN$ if $\p v_i$ lies in
the one-dimensional cone defined by $\tilde v_j$ and $r_i^j=0$ otherwise.
Our base space is the multiply weighted space determined by
\beq (\tilde z_1,\cdots,\tilde z_\Kt)\sim (\l^{\tilde w^1_j}
\tilde z_1,\cdots,\l^{\tilde w^\Kt_j}\tilde z_\Kt),~~~~~~
j=1,\cdots,\tilde K-\tilde n,  \eeql{ber} 
where $\tilde n=n-n\fib$ and the $\tilde w^i_j$ are any integers such that 
$\sum_i\tilde w^i_j \tilde v_i=0$.
The projection map from $V_\S$ (and, as we will see, from the Calabi--Yau
hypersurface) to the base is given by
\beq \tilde z_i=\prod_jz_j^{r_j^i}. \eeq
This is well defined: $z_j\to \l^{w_k^j}z_j$ leads to 
$\tilde z_i\to \l^{w_k^j r_j^i}\tilde z_i$ which is among the good
equivalence relations because applying $\p$ to $\sum w^j_k v_j=0$
gives $\sum w^j_k r_j^i \tilde v_i=0$.

A generic point in the base space will have $\tilde z_i\ne 0$ for all $i$,
implying $z_i\ne 0$ for all $v_i\not\in \D^*\fib$.
The choice of a specific point in $V_{\S\bas}$ and the use of all 
equivalence relations except for those involving only $v_i\in \D^*\fib$ 
allows to fix all $z_i$ except for those corresponding to 
$v_i\in \D^*\fib$.
Thus the preimage of a generic point in $V_{\S\bas}$
is indeed a variety in the moduli space determined by $\D^*\fib$.

What we have seen so far is just that $V_\S$ is a fibration over
$V_{\S\bas}$ with generic fiber $V_{\S\fib}$ (this is actually
the statement of an exercise on p. 41 of ref. \cite{Ful}) and
how this fibration structure manifests itself in terms of homogeneous
coordinates.
Now we also want to see how this can be extended to hypersurfaces.
To this end note that if $v_k\in\D^*\fib$ then $\<v_k,x\>$ only
depends on the equivalence class $[x]\in M\fib$ of $x$ under 
\beq x\sim y ~~~\hbox{ if }~~~ x-y \in M\bas . \eeq
Thus we may rewrite eq. (\ref{polyn}) as
\beq p=\sum_{[x]\in \D\fib\cap M\fib}a'_{[x]}
      \prod_{v_k\in\D^*\fib} z_k^{\<v_k,[x]\>+1}~~~\hbox{ with }~~~
    a'_{[x]}=\sum_{x\in [x]}a_x\prod_{v_k\not\in\D^*\fib}z_k^{\<v_k,x\>+1}.  
                                                                \eeql{polynf}
In each coordinate patch for $V_{\S\bas}$ this is just an equation
for the fiber with coefficients that are polynomial functions of
coordinates of the base space.

\del
There are two different occasions upon which the fiber may degenerate:
The fiber being a hypersurface in $V_{\S\fib}$, it can either happen
that $V_{\S\fib}$ itself degenerates or that the coefficients
of the equation determining the fiber hit a singular point in the moduli space.
While we do not know any way to read off the occurrence of the second case
from the toric data, we can surely see the first case: 
\enddel
Whenever a one-dimensional cone (with primitive generator $\tilde v_i$) 
in $\S\bas$ is the image of more than one one-dimensional cone in $\S$, 
the fiber becomes reducible over the divisor $\tilde z_i=0$ determined
by $v_i$.
Different components of the fiber correspond to different equations
$z_j=0$ with $\p v_j=r^i_j \tilde v_i$.
The intersection patterns of the different components of the
reducible fibers are crucial for understanding enhanced gauge symmetries
in type IIA string theory \cite{wit95,Arev} and F-theory \cite{Fv,MVI,MVII}.
Blowing down the corresponding subvarieties (and hence making the Calabi--Yau
space itself singular) leads to the appearance of 
non-perturbative enhanced gauge groups whose Lie algebras are determined
by the intersection patterns of the components of the fibers.
In terms of the $N$ lattice, the occurrence of enhanced gauge groups can
be easily inferred by studying the preimage of a one-dimensional cone in
$\S\bas$.
In particular, as noted by Candelas and Font \cite{cafo}, under
favorable circumstances the Dynkin diagram of the corresponding Lie
algebra can be seen directly in the toric diagram in the $N$ lattice.

\subsection{Fibrations and weight systems}


\def\fib{_{\rm f}} 

As for polyhedra, weight systems provide a very useful and economic tool
for constructing and describing fibrations. 
We only consider toric CY fibrations
which require a reflexive section through the origin of 
$\D^*\subset N_\IR$ whose dimension is equal to 
$n\fib$ for $(n\fib-1)$-dimensional CY-fibers. 
In the $M$ lattice this corresponds to a 
projection onto the dual reflexive polyhedron along an 
$(n-n\fib)$-dimensional subspace. 
Hence, on either side, we need to specify a linear subspace 
$(N\fib)_\IR\subset N_\IR$ or $(M\bas)_\IR\subset M_\IR$. 
This can be done, for example, by singling out 
vectors that span the subspace or by representing the subspace
as an intersection of hyperplanes.

If the polyhedron is given in terms of some CWS it is natural to try to
specify this linear subspace by using some part of the weight information. 
Note that a lower 
index $i$ in a CWS $n_i^{(j)}$ corresponds to a vertex of the minimal polytope
$\nabla$, or, by duality, to a bounding hyperplane in the $M$ lattice 
(which, if $\D$ is embedded into $\IR^k$, is the intersection of the 
coordinate hyperplane $a_i=0$, or $x_i=-1$, with the
affine or linear subspace that supports $\D$). 
An upper index $j$ corresponds to a
simplex in the CWS and, hence, to the linear subspace spanned by that simplex.
Actually this can be regarded as a special case of the former correspondence,
since the subspace that is spanned by a simplex $S^j$ is generated by those 
vectors of $\nabla$ for which $n_i^{(j)}\neq0$.

The simplest case is therefore the situation where a subset of the IP simplex
structure of the minimal polytope $\nabla$ provides a weight system for
$\D\fib$. 
In turn, we can engineer fibrations where the 
fiber corresponds to a certain weight system if we start by generating 
combined weight systems $\q$ with a given subsystem $\q\fib$. 
As usual, one has to check that $\D(\q)$ is reflexive, which in up 
to 4 dimensions is equivalent to the IP property 
(in the case of Calabi--Yau 4-fold CWS
with $\D(\q)$ IP but not reflexive we can proceed with 
reflexive subpolytopes $\D'\subset\D(\q)$). Because of the additional
equations that come from the extended CWS $\q$,
the set of solutions to eq. (\ref{AffineEq}) is smaller than those for 
$\q\fib$ when we disregard the additional coordinates. 
Hence the intersection of $\D^*$ with the subspace 
spanned by the vertices of $\nabla\fib$ contains 
$\D(\q\fib)^*$, but may be larger. We therefore obtain a fibration 
if the resulting $n\fib$-dimensional polytope 
is reflexive.%
\footnote{~
        If this is not the case we could proceed by dropping vertices of $\D$
        and trying to find a larger reflexive section, but this soon becomes 
        ugly and in view of the abundance of weight systems it is 
        hardly worth the effort.
}

In the case of elliptic fibrations this is always true: 
We want to show that $\D\fib^*=\D^* \cap (N\fib)_\IR$ is reflexive.
As we saw in section 4.1, $\D\fib$ can be identified with the image of 
$\D$ under $M_\IR\to (M\fib)_\IR$.
As the image of $\ipo_M$ under $M\to M\fib$ is an IP of $\D\fib$, $\D\fib$ 
has at least one IP.
Because of the well known fact that a two dimensional lattice polytope
is reflexive if it has precisely one IP, all that is left to show is that
$\D\fib$ cannot have more than one IP.
As $\D\fib^* \supseteq \nabla\fib$ implies $\D\fib \subseteq \nabla\fib^*$ 
and $\nabla\fib$ as a lattice polytope with a single IP is reflexive, 
$\nabla\fib^*$ and hence $\D\fib$ has precisely one IP, i.e. $\D\fib$ is
indeed reflexive.
\del
In the case of elliptic fibrations this is always true: 
Consider $\D\fib=(\D^* \cap (N\fib)_\IR)^*$, the image of $\D$ under 
$M_\IR\to (M\fib)_\IR$.
As the image of $\ipo$ is an IP of $\D\fib$, $\D\fib$ has at least one IP.
And $\D^* \cap (N\fib)_\IR \supseteq \nabla\fib$ implies 
$\D\fib \subseteq \nabla\fib^*$, 
but $\nabla\fib$ as a 2d IP polytope is reflexive, so $\nabla\fib^*$ and 
hence $\D\fib$ has precisely one IP, i.e. is reflexive. 
\enddel
\del
In the case of elliptic fibrations the last point is trivial because in
two dimensions all lattice polytopes with one interior point are reflexive.
Moreover, the intersection $\D^*\cap N\fib$ is always a lattice polytope,
which can be seen as follows: If $\D^*\cap N\fib$ contains an IP lattice 
polytope reflexivity of $\D^*$ implies that the 
intersection of any of its facets with $N\fib$ either is empty or is a 
hyperplane at distance one from the origin in $N\fib$. In 2 dimensions,
these intersections are the bounding lines of a polygon containing 
$\nabla\fib$. It is easy to see that all vertices
of such a polygon must be lattice points. 
{\rm\scriptsize(laesst sich, wenn ich mich nicht vertan habe, auch auf 
        K3-Faserungen erweitern, da dann die Facetten Polygone und das ganze
        noch immer uebersichtlich ist; schreib ich aber nicht rein, weil man 
        bei K3 sowieso den refcheck des Schnittes braucht)}
\enddel
To obtain elliptic fibrations in Weierstrass form, as they are mostly used
in F theory compactifications \cite{Fv,MVI,MVII}, we 
thus only need to take
the weight system $\q\fib=(1,2,3)/6$ as the fiber part of a
combined weight system and check for reflexivity of $\D(\q)$.

It is no surprise that the situation we described does not cover
the general case: 
Given an IP simplex structure of a minimal polytope 
$\nabla\fib\subseteq \D^*\fib$ 
it is not always possible to extend it to a simplex 
structure for a minimal polytope $\nabla\subseteq \D^*$. 
As a simple example in 4 dimensions 
we consider the points $V_1,\ldots,V_6\in N$ with coordinates
given by the columns of the matrix
\beq
        (V_1,\ldots,V_6)=\pmatrix{
        \;\;-1~~&~~1~~&1&-1~~&-1~~&-1~~\cr~~0 & 1 &-1~~& 0 & 0 & 0\cr
                ~~0 & 0 & 0 & 1 &-1~~&-1~~\cr   ~~0 & 0 & 0 & 0 & 1 &-1~~\cr}.
\eeq
The first three points provide a weight system (2,1,1) for an elliptic fiber
such that $\D^*\fib$ is supported by the 1--2 plane. 
$V_1$ is contained in the convex hull of $V_2,\ldots,V_6$, which
is a minimal polytope $\nabla$ as defined in section 2.1. 
The weight system for 
$\nabla$ is  $\q=(2,2,2,1,1)/8$ and we cannot use $(2,1,1)/4$
as part of a CWS corresponding to a minimal polyhedron $\nabla$ with an
IP simplex structure in the sense of section 2.
In this situation it makes sense to use generalized IP simplex structures
where the vertices of the IP simplices are lattice points (but not necessarily
vertices) of $\D^*$ and we do not insist on the non-redundancy implied by 
our original definition of an IP simplex structure.

Having made this point we may use the linear relations among the
vertices of the simplices $(V_1,V_2,V_3)$ and $(V_2,V_3,V_4,V_5,V_6)$ to 
arrive at the CWS $\q^{(1)}=(2,1,1,0,0,0)/4$ and 
$\q^{(2)}=(0,2,2,2,1,1)/8$. 
A CWS of this type was not considered in our classification scheme because 
$V_1=(V_2+V_3)/2$ is 
redundant when combined with the vertices that 
correspond to $\q^{(2)}$. It does, however, lead to a perfectly sensible
system of equations (\ref{AffineEq}), the convex hull of whose solutions 
is $\nabla^*$ (in our example all polytopes are simplices).
Actually, for $\q=(2,2,2,1,1)/8$ we find that $(\D(\q))^*$ has 
the seven lattice points $V_2,\ldots,V_6$, $\mathbf 0$ and $(-1,0,-1,0)^T$, 
but no reflexive subpolytope. The CWS 
$\{\q^{(1)},\q^{(2)}\}$, on the other
hand, leads to a polytope $(\D(\q^1,\q^2))^*$ with nine lattice points 
and the reflexive subpolytope that we started with: 
The addition of $V_1$ refines the lattice generated by the vertices of 
$\nabla$ in such a way that 
the convex hull of $V_2,\ldots,V_6$ on the finer $N$ lattice now contains the
additional lattice points $V_1$ and $-V_1$.
As an aside we thus observe that a CWS corresponding to a generalized IP
simplex structure may also be used
to encode certain sublattices of $M_{\rm finest}$.
Probably most polytopes can be directly specified by 
using a generalized IP simplex structure and the corresponding CWS.
A counterexample is given by the $\IZ_5$ quotient of the quintic at the end 
of section 3, where the $N$ lattice is not 
        generated by $\D^*\cap N$.
        But in practice such a representation is only useful if the number of 
        equations is small.
In any case combined weight systems provide a simple 
construction for toric fibrations and can always be used to specify reflexive
sections.


\subsection{Toric fibrations and weighted projective spaces}

There is a more subtle way in which a fibration structure can be encoded in
a weight system. It only works for codimension 1 fibers, but it is
quite interesting for historical and practical reasons.
When string dualities led to interest in K3 fibrations, the first
examples were constructed in the context of weighted projective spaces
\cite{klm,hly}. It turned out that
these examples are indeed special cases of toric fibrations in the sense
that they correspond to reflexive projections of Newton polyhedra of
transversal hypersurfaces in weighted $\IP^3$ or $\IP^4$.

Actually, reflexive objects of codimension 1 were first observed on these
Newton polyhedra, either as reflexive facets or as 
reflexive sections through the
IP in the $M$ lattice \cite{cafo}. Since what we really need for a toric 
fibration is a reflexive section in $N_\IR$, the question arises 
whether there is a reflexive projection of $\D$ onto one of 
its facets. A simple necessary
condition for this is provided by the following observations.
We work with the embedding space of lemma 4
and do not distinguish between objects in $M_\IR$ and their images under
the embedding map.
\\[3pt]
{\bf Lemma 10:}\\
a) For a polytope $\D$ defined by a weight
system $n_i$ only facets that are supported by a lattice hyperplane $x_l=-1$ 
can have interior points.\\
b) If ${\mathbf y}=(y_1,\ldots,y_k)\in M$ has $y_l=-1$, $y_i\ge 0$ for
$i\ne l$, the map $\p_{\mathbf y}: M_\IR\to M_\IR$,
$\p_{\mathbf y} {\mathbf x} = {\mathbf x} + (x_l+1) {\mathbf y}$ has the 
following properties:\\
It is a projection to the affine subspace $x_l=-1$, i.e. 
$\p_{\mathbf y}^2=\p_{\mathbf y}$ and 
$\p_{\mathbf y} M_\IR=M_\IR\cap \{x_l=-1\}$.\\
It respects the lattice structure, i.e. if ${\mathbf x}\in M$ then 
$\p_{\mathbf y} {\mathbf x}\in M\cap \{x_l=-1\}$,\\
The image of $\D$ is the corresponding facet of $\D$, i.e. 
$\p_{\mathbf y}\D=\D\cap \{x_l=-1\}$.\\
c) There is a one-to-one correspondence between maps $\p$ with the same
properties as in b) such that $\ipo$ gets mapped to an IP of the facet 
with $x_l=-1$ and partitions of the weight $n_l$ by the remaining 
weights, i.e. 
$n_l=\sum_{i\neq l} y_in_i$ where the 
$y_i$ are nonnegative integers.
\\[2pt]
{\it Proof:} a) If an interior point of a facet is not on some hyperplane 
$x_l=-1$ all $x_i$ must be nonnegative, but this is only possible for the 
interior point of $\D$.\\
b) The first two statements follow directly from the definition of 
$\p_{\mathbf y}$. $\p_{\mathbf y}\D\supseteq\D\cap \{x_l=-1\}$ follows from
the fact that $\p_{\mathbf y}$ is a projection. For 
$\p_{\mathbf y}\D\subseteq\D\cap \{x_l=-1\}$  we note that $\D\cap \{x_l=-1\}$
is the convex hull of the lattice points in $\{x_l=-1\}$ with $x_i\ge -1$
and that every vertex of $\D$ gets mapped to such a lattice point.\\
c) If $\p$ is such a map, then we choose $\mathbf y$ to be the IP of the facet 
to which $\ipo$ is mapped (this implies $y_i\ge 0$ for
$i\ne l$), and the partition of $n_l$ follows from the fact
that $\sum y_in_i=0$. Conversely, if $\mathbf y$ is defined by such a
partition we have to show that it is interior to the facet. This follows from
the facts that $\ipo$ is interior to $\D$ and $\mathbf y=\p_{\mathbf y}\ipo$.
 \hBo
\del
b) follows from eq. (\ref{AffineEq}), and c) follows from the simple form 
of the projection
vector, which guarantees that the projected lattice points of $\D$ have 
integral and nonnegative coordinates $a_i$.
\enddel

A necessary condition for the existence of a reflexive projection of $\D$
onto one of its facets is therefore that one of the weights $n_i$ has
a {\it unique} partition in terms of the other weights. Using this criterion 
we found all such projections for single weight systems with $k\le5$ by first
searching for weights with unique partitions and then checking reflexivity
of the corresponding facets.
The results are given in table 3 for the case of elliptic K3 surfaces and
they are available on our web page \cite{KScy} for K3-fibered Calabi--Yau
manifolds (cf. table 4). 

We can find a set of generators for $(N\fib)_\IR$
by solving the equation $\<V,{\mathbf y}\>=0$ for a general
linear combination $V=\sum c_j V_j$ of the vertices $V_j$ of $\nabla$. 
Since $\<V_i,{\mathbf y}\>=y_i$ we obtain the solutions $V'_i=V_i+y_iV_l$ for 
$i\neq l$.
The linear relations among the $V_i'$ are given by the corresponding subset of
the original weights. 
In general 
they do not provide a weight system for the fiber because
the points $V_i'$ need not belong to $\D^*$.

This is easy to see for the class of weights $(1,1,2n_3,2n_4,2n_5)$ 
that was considered by Klemm, Lerche and Mayr \cite{klm}. Here $V_2'=V_1+V_2$
and $\<V_2',{\mathbf x}\>=x_1+x_2$ for 
${\mathbf x}\in M$ with coordinates $x_i$.
But $\sum x_in_i=0$ implies that $x_1+x_2$
is even, so $V_2'$ is not a primitive lattice vector in $N$ and can be
divided by 2, which leads to the weight system $(1,n_3,n_4,n_5)$ for the K3
fiber. The slightly more complicated example $(8,4,3,27,42)/84$ was 
given by Hosono, Lian and Yau \cite{hly}. The first weight has a unique 
partition with $y_2=2$ and $y_3=y_4=y_5=0$, so that $(N\fib)_\IR$ 
is spanned by $V'_2=V_2+2V_1$ and $V_i$ with $i>2$. This time 
$8x_1+4x_2+3x_3 +27x_4+42x_5=0$ implies that $x_2+2x_1$ is a multiple of 3
and the primitive lattice vector $V'_2/3\in\D^*$ leads to the weight system
$(4,1,9,14)$ for the fiber, which agrees with the normalized weights for the
fiber given in \cite{hly}.
If more of the coefficients $y_i$ are nonvanishing, it is, of course, still
possible to compute a weight system for the fiber, but this gets more
tedious and we would also lose the direct connection with the original 
weights or we would have to introduce many redundant coordinates in a CWS.

Another strategy for identifying reflexive projections of $\D(q)$ that can
be used in the codimension 1 case follows from the fact that such a projection
either must be along a line parallel to a facet or onto that facet whenever
a facet has an interior point. If the number of facets with interior points
is large enough this allows us to find all reflexive projections. The result
of this analysis is indicated in the next-to-last column of 
tables 3 and 4. 
        The K3 surfaces in table 2 are all elliptic, since their combinded
        weight systems contain 2 dimensional subsystems. 

With our strategies to identify reflexive projections (onto facets) we 
generalized the results of \cite{klm,hly} and
extended the scope from the transversal case to the complete
list of 184026 IP weight systems, where we could identify 124701 fibrations.
The efficiency of our approach can be inferred from the fact that
we found 5370 fibrations 
for the 7555 transversal cases,
wheras only 628 fibrations yielded to the methods of \cite{hly}.

\bigskip

{\it Acknowledgements:} M.K. is partly supported by the Austrian 
Research Funds FWF grant Nr. P11582-PHY. 
The research of H.S. is supported by the European Union 
TMR project ERBFMRX-CT-96-0045.

\newpage

\section*{Appendix: Various tables}

\BC{\small\ni
\begin{tabular}{|l|r|r|r|r|} \hline\hline
IP simplex structure & total  & span & lp-min.  & r-min. \\ \hline\hline
$S_1=V_1V_2V_3V_4V_5$ & 184026 & 38730 & 16437 & 206\\
$S_1=V_1V_2V_3V_4,\; S_2=V_1V_2V_3'V_4'$ & 16040 & 6365 & 143 & 51\\
$S_1=V_1V_2V_3V_4,\; S_2=V_1V_2'V_3'$ & 1122 & 727 & 40 & 29\\
$S_1=V_1V_2V_3,\; S_2=V_1'V_2'V_3'$ & 6 & 6  & 3 & 3\\
$S_1=V_1V_2V_3,\; S_2=V_1V_2'V_3',\; S_3=V_1V_2''V_3''$ & 36 & 36 & 4 & 4\\
$S_1,\ldots,S_{m-1}\hbox{ as for }n=3,\;S_m=V_1^{(m)}V_2^{(m)}$ & 
                                               116 & 79 & 19 & 15\\ \hline
total & 201346 & 45943 & 16646 & 308 \\ \hline\hline 
\end{tabular}}
\\[12pt]{\bf Table 1}: 
IP simplex structures and numbers of corresponding IP CWS for $n=4$
\\[-2pt](for $n=3$, see lemma 3).\HS190 ~
\EC

\noindent
\begin{center}  {\footnotesize \def\NWS{&&&&\\[-5.5pt]}\def\NCWS{\\[-1pt]}
\HS-15 \begin{tabular}{|r|rrrrr|c|cc|cc|} \hline\hline
d & $n_1$  & $n_2$  & $n_3$  & $n_4$  & $n_5$  & & 
                        $P$ & $V$ & $\5 P$ &$\5 V$    \\ \hline\hline
3 & 1 & 1 & 1 & 0 & 0&\NWS   3 & 1 & 0 & 0 & 1 & 1 & r&30 &5 &6 &5\NCWS \hline
3 & 1 & 1 & 1 & 0 & 0&\NWS   4 & 2 & 0 & 0 & 1 & 1 & r&31 &6 &7 &5\NCWS \hline
3 & 1 & 1 & 1 & 0 & 0&\NWS   4 & 1 & 0 & 0 & 2 & 1 & s&23 &7 &8 &6 \NCWS \hline
3 & 1 & 1 & 1 & 0 & 0&\NWS   6 & 3 & 0 & 0 & 2 & 1 & s&24 &6 &9 &5 \NCWS \hline
3 & 1 & 1 & 1 & 0 & 0&\NWS   6 & 2 & 0 & 0 & 3 & 1 & s&21 &5 &9 &5 \NCWS \hline
3 & 1 & 1 & 1 & 0 & 0&\NWS   6 & 1 & 0 & 0 & 3 & 2 & s&14 &7 &11 &6 \NCWS\hline
4 & 2 & 1 & 1 & 0 & 0&\NWS   4 & 2 & 0 & 0 & 1 & 1 & r&35 &5 &7 &5 \NCWS \hline
4 & 2 & 1 & 1 & 0 & 0&\NWS   4 & 1 & 0 & 0 & 2 & 1 & s&23 &6 &9 &5 \NCWS \hline
4 & 2 & 1 & 1 & 0 & 0&\NWS   6 & 3 & 0 & 0 & 2 & 1 & s&27 &5 &9 &5 \NCWS \hline
4 & 1 & 2 & 1 & 0 & 0&\NWS   4 & 1 & 0 & 0 & 2 & 1 & s&19 &5 &9 &5 \NCWS \hline
4 & 1 & 2 & 1 & 0 & 0&\NWS   6 & 3 & 0 & 0 & 2 & 1 & s&18 &6 &12 &5 \NCWS\hline
\hline \end{tabular}
~~~~~~
\begin{tabular}{|r|rrrrrr|c|cc|cc|} \hline\hline
d & $n_1$  & $n_2$  & $n_3$  & $n_4$  & $n_5$  & $n_6$  & & 
                        $P$ & $V$ & $\5 P$ &$\5 V$    \\ \hline\hline
4 & 1 & 2 & 1 & 0 & 0&&\NWS   6 & 2 & 0 & 0 & 3 & 1 &&s&16 &6 &14 &6\NCWS\hline
4 & 1 & 2 & 1 & 0 & 0&&\NWS   6 & 1 & 0 & 0 & 3 & 2 &&s&12 &6 &14 &6\NCWS\hline
6 & 3 & 2 & 1 & 0 & 0&&\NWS   6 & 3 & 0 & 0 & 2 & 1 &&s&21 &5 &12 &5\NCWS\hline
6 & 2 & 3 & 1 & 0 & 0&&\NWS   6 & 2 & 0 & 0 & 3 & 1 &&s&15 &5 &15 &5\NCWS\hline
6 & 2 & 3 & 1 & 0 & 0&&\NWS   6 & 1 & 0 & 0 & 3 & 2 &&s&10 &6 &20 &6\NCWS\hline
6 & 1 & 3 & 2 & 0 & 0&&\NWS   6 & 1 & 0 & 0 & 3 & 2 &&s&9 &5 &18 &5\NCWS\hline
\hline
3 & 1 & 1 & 1 & 0 & 0&&\NWS   2 & 0 & 0 & 0 & 1 & 1 && r&30 &6 &6 &5\NCWS\hline
4 & 2 & 1 & 1 & 0 & 0&&\NWS   2 & 0 & 0 & 0 & 1 & 1 && r&27 &6 &7 &5\NCWS\hline
6 & 3 & 2 & 1 & 0 & 0&&\NWS   2 & 0 & 0 & 0 & 1 & 1 && s&21 &6 &9 &5\NCWS\hline
\hline
 &  &  &  &  &  &  & &&&&\\[-9pt]
2 & 1 & 1 & 0 & 0 & 0 & 0&\NWS   
2 & 0 & 0 & 1 & 1 & 0 & 0& r &27 &8 &7 &6\\[-5.5pt] 
2 & 0 & 0 & 0 & 0 & 1 & 1 &&&&&\\[3pt] \hline
\hline 
\end{tabular}\HS-19 ~}
\\[12pt]{\bf Table 2:} IP CWS for $n=3$. 
        The columns indicate the minimality type (`s' for span, `l' for 
\\\HS22 lp-minimality and `r' for r-minimality) and point and vertex numbers 
        for $\D$ and $\D^*$. As 
\\\HS12 r-minimality implies lp-minimality and the latter implies the span 
        property for $n=3$, 
\\      we have given only the strongest statement in each case.\HS139 ~
\end{center}

\newpage

\begin{center}  { \footnotesize   \def\NLHL{\\[-5pt]} \def\N{--}

\def\tabtop
{\begin{tabular}{||\TVR{3.5}1 
                c|c@{~}c@{~}c@{~}c|c|r@{~}c@{~}|r@{~}c@{~}|c@{~~~}c||} 
        \hline\hline    \VR41
$d$     &$n_1$  & $n_2$ & $n_3$ & $n_4$ &  
        & ~$P$ & ~~$V$~ & ~$\5P$ & ~~$\5V$~ & $\P$      
        & F \\\hline\hline
}                                                       \hspace*{-6mm}\tabtop
4 & 1 &1 &1 &1 &        r& 35 & 4 & 5 & 4 & 0 & 0\NLHL
5 & 1 &1 &1 &2 &        r& 34 & 6 & 6 & 5 & 0 & 0\NLHL
6 & 1 &1 &1 &3 &        r& 39 & 4 & 6 & 4 & 0 & 0\NLHL
6 & \bf1 &1 &2 &2 &     r& 30 & 4 & 6 & 4 & 1 & 1\NLHL
7 & \bf1 &1 &2 &3 &     r& 31 & 7 & 8 & 6 & 1 & 1\NLHL
8 & 1 &2 &2 &3 &        s& 24 & 6 & 8 & 5 & 0 & 0\NLHL
8 & \bf1 &1 &2 &4 &     r& 35 & 4 & 7 & 4 & 1 & 1\NLHL
9 & 1 &\bf2 &3 &3 &     s& 23 & 6 & 8 & 5 & 1 & 1\NLHL
9 & \bf1 &1 &3 &4 &     r& 33 & 5 & 9 & 5 & 1 & 1\NLHL
10 & 1 &2 &2 &5 &       s& 28 & 4 & 8 & 4 & 0 & 0\NLHL
10 & 1 &\bf2 &3 &4 &    s& 23 & 7 & 11 & 6 & 1 & 1\NLHL
10 & \bf1 &1 &3 &5 &    r& 36 & 5 & 9 & 5 & 1 & 1\NLHL
11 & 1 &\bf2 &3 &5 &    l& 24 & 8 & 13 & 7 & 1 & 1\NLHL
12 & 1 &\bf2 &3 &6 &    s& 27 & 4 & 9 & 4 & 1 & 1\NLHL
12 & 1 &\bf2 &4 &5 &    s& 24 & 5 & 12 & 5 & 1 & 1\NLHL
12 & 1 &\bf3 &4 &4 &    s& 21 & 4 & 9 & 4 & 1 & 1\NLHL
12 & 2 &\bf3 &3 &\bf4 & s& 15 & 4 & 9 & 4 & 2 & 2\NLHL
12 & \bf1 &1 &4 &6 &    r& 39 & 4 & 9 & 4 & 1 & 1\NLHL
12 & \bf2 &2 &3 &5 &    s& 17 & 5 & 11 & 5 & 1 & 1\NLHL
13 & 1 &\bf3 &4 &5 &    l& 20 & 7 & 15 & 7 & 1 & 1\NLHL
14 & 1 &\bf2 &4 &7 &    s& 27 & 5 & 12 & 5 & 1 & 1\NLHL
14 & 2 &3 &\bf4 &\bf5 & s& 13 & 7 & 16 & 7 & 3 & 2\NLHL
14 & \bf2 &2 &3 &7 &    s& 19 & 5 & 11 & 5 & 1 & 1\NLHL
15 & 1 &\bf2 &5 &7 &    s& 26 & 6 & 17 & 6 & 1 & 1\NLHL
15 & 1 &\bf3 &4 &7 &    s& 22 & 6 & 17 & 6 & 1 & 1\NLHL
15 & 1 &\bf3 &5 &6 &    s& 21 & 5 & 15 & 5 & 1 & 1\NLHL
15 & 2 &3 &5 &5 &       s& 14 & 6 & 11 & 5 & 1 & 0\NLHL
15 & \bf3 &3 &4 &5 &    s& 12 & 5 & 12 & 5 & 1 & 1\NLHL
16 & 1 &\bf2 &5 &8 &    s& 28 & 5 & 14 & 5 & 1 & 1\NLHL
16 & 1 &\bf3 &4 &8 &    s& 24 & 5 & 12 & 5 & 1 & 1\NLHL
16 & 1 &\bf4 &5 &6 &    s& 19 & 6 & 17 & 6 & 1 & 1\NLHL
16 & 2 &3 &\bf4 &7 &    s& 14 & 6 & 18 & 6 & 2 & 1\NLHL
17 & 2 &3 &\bf5 &7 &    l& 13 & 8 & 20 & 8 & 2 & 1\NLHL
18 & 1 &\bf2 &6 &9 &    s& 30 & 4 & 12 & 4 & 1 & 1\NLHL
18 & 1 &\bf3 &5 &9 &    s& 24 & 5 & 15 & 5 & 1 & 1\NLHL
18 & 1 &\bf4 &6 &7 &    s& 19 & 6 & 20 & 6 & 1 & 1\NLHL
18 & 2 &3 &\bf4 &9 &    s& 16 & 5 & 14 & 5 & 2 & 1\NLHL
18 & 2 &3 &\bf5 &8 &    s& 14 & 6 & 20 & 6 & 2 & 1\NLHL
18 & 3 &4 &5 &\bf6 &    s& 10 & 6 & 17 & 6 & ? & 1\NLHL
19 & 3 &4 &5 &\bf7 &    l& 9 & 7 & 24 & 8 & ? & 1\NLHL
20 & 1 &\bf4 &5 &10 &   s& 23 & 4 & 13 & 4 & 1 & 1\NLHL
20 & 2 &3 &\bf5 &10 &   s& 16 & 5 & 14 & 5 & 2 & 1\NLHL
20 & 2 &5 &\bf6 &\bf7 &\N& 11 & 5 & 23 & 5 & 3 & 2\NLHL
20 & 2 &\bf4 &5 &9 &\N& 13 & 4 & 23 & 4 & 1 & 1\NLHL
20 & 3 &4 &5 &8 &       s& 10 & 6 & 22 & 6 & ? & 0\NLHL
21 & 1 &\bf3 &7 &10 &\N& 24 & 4 & 24 & 4 & 1 & 1\NLHL
21 & 1 &\bf5 &7 &8 &\N& 18 & 5 & 24 & 5 & 1 & 1\NLHL
21 & 2 &3 &\bf7 &9 &    s& 14 & 6 & 23 & 6 & 2 & 1\\
\hline\end{tabular}     \hspace{9mm}    
%
        \tabtop
21 & 3 &5 &\bf6 &7 &\N& 9 & 5 & 21 & 5 & ? & 1\NLHL
22 & 1 &\bf3 &7 &11 &\N& 25 & 5 & 20 & 5 & 1 & 1\NLHL
22 & 1 &\bf4 &6 &11 &   s& 22 & 6 & 20 & 6 & 1 & 1\NLHL
22 & 2 &\bf4 &5 &11 &\N& 14 & 5 & 19 & 5 & 1 & 1\NLHL
24 & 1 &\bf3 &8 &12 &   s& 27 & 4 & 15 & 4 & 1 & 1\NLHL
24 & 1 &\bf6 &8 &9 &    s& 18 & 5 & 24 & 5 & 1 & 1\NLHL
24 & 2 &3 &8 &11 &\N& 15 & 4 & 27 & 4 & 1 & 0\NLHL
24 & 2 &3 &\bf7 &12 &   s& 16 & 5 & 20 & 5 & 2 & 1\NLHL
24 & 3 &4 &5 &12 &      s& 12 & 5 & 18 & 5 & ? & 0\NLHL
24 & 3 &4 &\bf7 &10 &   s& 10 & 5 & 26 & 6 & 2 & 1\NLHL
24 & 3 &\bf6 &7 &8 &\N& 9 & 4 & 21 & 4 & ? & 1\NLHL
24 & 4 &5 &6 &\bf9 & s& 8 & 5 & 26 & 6 & ? & 1\NLHL
25 & 4 &5 &7 &\bf9 &\N& 7 & 5 & 32 & 6 & ? & 1\NLHL
26 & 1 &\bf5 &7 &13 &\N& 21 & 5 & 24 & 5 & 1 & 1\NLHL
26 & 2 &3 &8 &13 &\N& 16 & 5 & 23 & 5 & 1 & 0\NLHL
26 & 2 &5 &\bf6 &13 &\N& 13 & 5 & 23 & 5 & 2 & 1\NLHL
27 & 2 &5 &\bf9 &11 &\N& 11 & 6 & 32 & 6 & 2 & 1\NLHL
27 & 5 &6 &7 &9 &\N& 6 & 5 & 30 & 6 & ? & 0\NLHL
28 & 1 &\bf4 &9 &14 &\N& 24 & 4 & 24 & 4 & 1 & 1\NLHL
28 & 3 &4 &\bf7 &14 &   s& 12 & 5 & 18 & 5 & 2 & 1\NLHL
28 & 4 &6 &7 &\bf11 &\N& 7 & 4 & 35 & 4 & ? & 1\NLHL
30 & 1 &\bf4 &10 &15 &  s& 25 & 5 & 20 & 5 & 1 & 1\NLHL
30 & 1 &\bf6 &8 &15 &   s& 21 & 5 & 24 & 5 & 1 & 1\NLHL
30 & 2 &3 &10 &15 & s& 18 & 4 & 18 & 4 & 1 & 0\NLHL
30 & 2 &\bf6 &7 &15 &\N& 13 & 4 & 23 & 4 & 1 & 1\NLHL
30 & 3 &4 &\bf10 &13 &\N& 10 & 5 & 35 & 5 & 2 & 1\NLHL
30 & 4 &5 &6 &15 &      s& 10 & 5 & 20 & 5 & ? & 0\NLHL
30 & 5 &6 &8 &\bf11 &\N& 6 & 4 & 39 & 4 & ? & 1\NLHL
32 & 2 &5 &\bf9 &16 &\N& 13 & 5 & 29 & 5 & 2 & 1\NLHL
32 & 4 &5 &7 &16 &\N& 9 & 5 & 27 & 5 & ? & 0\NLHL
33 & 3 &5 &\bf11 &14 &\N& 9 & 4 & 39 & 4 & 2 & 1\NLHL
34 & 3 &4 &\bf10 &17 &\N& 11 & 6 & 31 & 6 & 2 & 1\NLHL
34 & 4 &6 &7 &\bf17 &\N& 8 & 5 & 31 & 5 & ? & 1\NLHL
36 & 1 &\bf5 &12 &18 &\N& 24 & 4 & 24 & 4 & 1 & 1\NLHL
36 & 3 &4 &\bf11 &18 &\N& 12 & 4 & 30 & 4 & 2 & 1\NLHL
36 & 7 &8 &9 &12 &\N& 5 & 4 & 35 & 4 & ? & 0\NLHL
38 & 3 &5 &\bf11 &19 &\N& 10 & 5 & 35 & 5 & 2 & 1\NLHL
38 & 5 &6 &8 &\bf19 &\N& 7 & 5 & 35 & 5 & ? & 1\NLHL
40 & 5 &7 &8 &20 &\N& 8 & 4 & 28 & 4 & ? & 0\NLHL
42 & 1 &\bf6 &14 &21 &  s& 24 & 4 & 24 & 4 & 1 & 1\NLHL
42 & 2 &5 &14 &21 &\N& 15 & 4 & 27 & 4 & 1 & 0\NLHL
42 & 3 &4 &\bf14 &21 &  s& 13 & 5 & 26 & 5 & 2 & 1\NLHL
44 & 4 &5 &\bf13 &22 &\N& 9 & 4 & 39 & 4 & 2 & 1\NLHL
48 & 3 &5 &\bf16 &24 &\N& 12 & 4 & 30 & 4 & 2 & 1\NLHL
50 & 7 &8 &10 &\bf25 &\N& 6 & 4 & 39 & 4 & ? & 1\NLHL
54 & 4 &5 &\bf18 &27 &\N& 10 & 5 & 35 & 5 & 2 & 1\NLHL
66 & 5 &6 &\bf22 &33 &\N& 9 & 4 & 39 & 4 & 2 & 1\NLHL
                                &&&&&&&&&&&     \\
\hline\end{tabular}\hspace*{-6mm}
}
\\[4mm]
\end{center}
        {\bf Table 3:}  The 95 K3 weight systems: $r,l,s$ denote the
                minimality type as in table 2, 
                $\P$ is the number of reflexive projections
                (if known) and $F$ denotes the number of reflexive
                projections onto facets. 
                The corresponding weights with unique
                partitions are indicated with bold face.

\newpage

\def\Rls{{\bf{--}}\,ls} \def\Rs{{\bf{--}}\,\,s} \def\Rr{{\bf{--}}\,\,r}
\def\Rl{{\bf{--}}\,\,l} \def\Rn{\bf{--}\,\bf{--}}  
\def\Tr{T\,r} \def\Ts{T\,s} \def\Tn{T\,\bf{--}}

\begin{center}          {\scriptsize 
\vspace*{-12mm}
\begin{tabular}{||\TVR{3.4}1 c|ccccc|l|cc|cc|cc|cc||} \hline\hline      \VR41
$d$     &$n_1$  & $n_2$ & $n_3$ & $n_4$ & $n_5$ &  TM   & $h_{11}$ & $h_{12}$ 
& $P$ & $V$ & $\5P$       & $\5V$       & $\P$  & F \\\hline\hline
47 &3 & 4 & 5 & 14 & 21 & \Rls & 26&39 & 54& 18 & 35&15 & ?&0\\\hline
69 &7 & 8 & 10 & 19 & 25 & \Rls& 59&10 & 16&13 & 75&21 & ?&0\\
97 &7 & 8 & 11 & \bf26 & 45 & \Rls& 63&15 & 24&15 & 71&21 & ?&1\\\hline
84 &\bf1 & 1 & 12 & 28 & 42 & \Tr& 11&491 & 680&5 & 26&5 & 1&1\\\hline
280 &7 & 19 & \bf40 & 87 & 127 & \Rn& 491&11 & 26&5 & 680&5 & 2&1\\\hline
24 &3 & 4 & 5 & 6 & 6 & \Ts& 10&34 & 36&8 & 12&7 & ?&0\\
26 &3 & 4 & 5 & 7 & 7 & \Rls& 22&22 & 31&13 & 21&10 & ?&0\\
33 &3 & 6 & 6 & 7 & 11 & \Rn& 19&37 & 34&7 & 22&6 & ?&0\\
36 &3 & 6 & 6 & 10 & 11 & \Tn& 19&49 & 38&7 & 22&6 & ?&0\\\hline
26 &3 & 4 & 5 & \bf6 & 8 & \Rs& 14&24 & 32&14 & 19&10 & ?&1\\
36 &5 & \bf7 & 7 & 8 & 9 & \Rn& 30&12 & 19&10 & 28&9 & ?&1\\
39 &3 & \bf6 & 9 & 10 & 11 & \Ts& 17&41 & 33&12 & 22&13 & ?&1\\
52 &4 & 6 & \bf8 & 11 & 23 & \Tn& 29&33 & 34&9 & 36&8 & ?&1\\\hline
34 &3 & \bf6 & 7 & 8 & \bf10 & \Rs& 18&20 & 27&13 & 23&12 & ?&2\\
44 &4 & \bf8 & 9 & 10 & \bf13 & \Rn& 29&17 & 22&9 & 31&9 & ?&2\\
55 &3 & 10 & \bf13 & 14 & \bf15 & \Ts& 28&16 & 23&12 & 35&14 & ?&2\\
63 &7 & 9 & \bf14 & 15 & \bf18 & \Tn& 44&8 & 15&6 & 37&6 & ?&2\\\hline
5 &1 & 1 & 1 & 1 & 1 & \Tr& 1&101 & 126&5 & 6&5 & 0&0\\
10 &1 & 1 & 1 & 3 & 4 & \Rr& 4&126 & 165&10 & 9&7 & 0&0\\
25 &1 & 5 & 5 & 6 & 8 & \Tn& 17&49 & 65&7 & 15&7 & 0&0\\
26 &1 & 5 & 5 & 7 & 8 & \Rl& 19&49 & 65&9 & 19&7 & 0&0\\\hline
20 &2 & 3 & 4 & 4 & 7 & \Rs& 13&45 & 51&10 & 14&8 & 1&0\\
20 &2 & 3 & 5 & 5 & 5 & \Ts& 6&48 & 50&8 & 11&6 & 1&0\\
30 &2 & 5 & 6 & 6 & 11 & \Rn& 27&39 & 45&7 & 25&6 & 1&0\\
36 &2 & 5 & 6 & 6 & 17 & \Tn& 24&54 & 60&7 & 25&6 & 1&0\\\hline
8 &\bf1 & 1 & 2 & 2 & 2 & \Tr& 2&86 & 105&5 & 7&5 & 1&1\\
13 &\bf1 & 1 & 2 & 4 & 5 & \Rr& 6&108 & 141&12 & 11&8 & 1&1\\
35 &2 & 7 & \bf8 & 9 & 9 & \Rn& 35&23 & 33&11 & 30&8 & 1&1\\
40 &4 & 5 & \bf9 & 10 & 12 & \Tn& 22&18 & 25&7 & 20&7 & 1&1\\\hline
19 &2 & 3 & \bf4 & 5 & 5 & \Rls& 11&33 & 43&14 & 14&9 & 2&1\\
27 &2 & 3 & \bf4 & 9 & 9 & \Ts& 14&44 & 56&9 & 13&7 & 2&1\\
36 &\bf4 & 4 & 6 & 9 & 13 & \Rn& 31&31 & 33&6 & 29&6 & 2&1\\
40 &\bf4 & 4 & 6 & 9 & 17 & \Tn& 26&38 & 39&7 & 29&6 & 2&1\\\hline
14 &\bf2 & 2 & \bf3 & 3 & 4 & \Tr& 5&51 & 57&10 & 10&7 & 2&2\\
19 &2 & \bf3 & 3 & \bf4 & 7 & \Rls& 11&39 & 51&14 & 16&9 & 2&2\\
30 &3 & \bf5 & 5 & \bf6 & 11 & \Rn& 33&21 & 33&6 & 25&6 & 2&2\\
35 &3 & \bf5 & 5 & \bf6 & 16 & \Tn& 26&28 & 42&7 & 25&6 & 2&2\\\hline
28 &4 & \bf5 & 5 & 6 & \bf8 & \Rs& 18&20 & 27&10 & 18&8 & 3&2\\
36 &4 & 6 & \bf8 & \bf9 & 9 & \Ts& 23&23 & 26&6 & 16&6 & 3&2\\
40 &4 & \bf7 & 7 & 10 & \bf12 & \Rn& 28&16 & 23&7 & 25&6 & 3&2\\
42 &6 & \bf7 & 7 & 10 & \bf12 & \Tn& 35&11 & 19&6 & 23&6 & 3&2\\
\hline\end{tabular}}
\\[5mm]
        {\bf Table 4:}  Examples from our list of 184026 IP weights 
        \cite{KScy} with various data including\\
~\HS33  Hodge numbers, point and vertex numbers, and numbers of reflexive 
        projections \\
~\HS15  (onto facets). $T$ indicates transversality and $M$ denotes the 
        minimality type.
\end{center}

\bye